\documentclass [12pt,oneside] {amsart}
\usepackage{latexsym,amsfonts,amssymb,euler,euscript,amscd,epsfig,isolatin1,
supertabular,rotating,epic,fullpage}

\theoremstyle{plain}
\newtheorem {Thm} {Theorem}[section]
\newtheorem* {Thm*} {Theorem}
\newtheorem* {Prop*} {Proposition}
\newtheorem {Lem}[Thm] {Lemma}
\newtheorem {Prop}[Thm] {Proposition}
\newtheorem {Cor}[Thm] {Corollary}

\theoremstyle {definition}
\newtheorem {Def}[Thm] {Definition}

\newtheorem {Rem}[Thm] {Remark}
\newtheorem {Exa}[Thm] {Example}
\newenvironment{Pf}[1]{{\noindent\sc Proof #1:}}{\qed\\}
\newcommand {\longto} {\longrightarrow}
\newcommand {\specialmap} [4] {\text {$ #1\negmedspace : #2 #3 #4 $}}
\newcommand {\map} [3] {\specialmap {#1} {#2}{\to} {#3}}
\newcommand {\longmap} [3] {\specialmap {#1} {#2}{\longto} {#3}}

\newcommand {\Aut}{\operatorname{Aut}}
\newcommand {\Hom} {\operatorname {Hom}}
\newcommand {\id} {\operatorname{id}}

\newcommand {\ind} {\operatorname{ind}}
\renewcommand {\(} {\left(}
\renewcommand {\)} {\right)} 

\renewcommand {\geq} {\geqslant}
\newcommand {\Oplus} {{\bigoplus\limits_{n\geq 0}}}
\newcommand {\cOplus} {\hat\Oplus}
\newcommand {\cOMU} {{\cOplus\MU_\Sn(X^n)t^n}}
\newcommand {\cOK} {{\cOplus K_\Sn(X^n)t^n}}

\def\acted{\hspace{.1cm}{
        \setlength{\unitlength}{.32mm}
        \linethickness{.09mm}
        \begin{picture}(8,8)(0,0)
        \qbezier(1,6)(3.5,8.3)(6,7)
        \qbezier(6,7)(9.5,4)(6,1)
        \qbezier(6,1)(3.5,-.3)(1,2)
        \qbezier(1,6)(1.9,7.5)(1.2,9)
        \qbezier(1,6)(3,6.1)(3.8,4.4)
        \end{picture}\hspace{.1cm}
        }}

\newcommand {\sub} {\subseteq}

\newcommand {\CC} {\mathbb C}
\def\Pfeil{{\unitlength 1em\begin{picture}(0,.1)
\put(0,.1){\vector(1,0){1.5}}\end{picture}}}
\def\colimname{{\unitlength.1em
\raisebox{-2.7\unitlength}{\begin{picture}(15.5,9.5)(0,0)
\put(0,2.7){$\operatorname{lim}$}
\put(.05,-.1){\Pfeil} 
\end{picture}}}}
\def\Colim{\mathop{\colimname}}
\newcommand {\on}[1] {\operatorname{#1}}
\newcommand {\BG} {{\on{BG}}}
\newcommand {\BGH} {{\on{B\mathcal G}}}

\newcommand {\Borel} {\on{EG}_+\wedge_G}
\newcommand {\borel}[1] {{\on{EG}\times_G #1}}

\newcommand {\EBG} {{E^0(\BG)}}
\newcommand {\EG}[1] {{E_G(#1)}}

\newcommand {\eps} {\varepsilon}

\newcommand {\eval} {{\on{eval}_\alpha}}

\newcommand {\forb}[1] {{\phi_{\on{orb}}(#1)}}
\renewcommand {\gg} {{(g_1,\dots,g_h)}}
\newcommand {\GG}[1] {{\frac{1}{|#1|}}}
\newcommand {\gpd} {{\mathcal G}}
\newcommand {\GS} {{G\wr\Sn}}
\newcommand {\hinf} {{H_\infty}}
\newcommand {\htpp} {$h$-tuple of commuting elements of $p$-power
        order }
\newcommand {\htpps} {$h$-tuples of commuting elements of $p$-power
        order }
\renewcommand {\ind}[2] {\on{ind}\arrowvert_{#1}^{#2}}

\newcommand {\HKR} {Hopkins-Kuhn-Ravenel }
\newcommand {\Kn} {{K(h)}}
\newcommand {\Korb}[1] {K_{\on{orb}}(#1)}

\newcommand {\LT}[1] {{E_h^0(#1)}}
\renewcommand {\ll} {\underline{l}}
\newcommand {\MG} {{M\acted G}}
\newcommand {\Mn} {{M^n\acted\Sn}}

\newcommand {\mmod} {/\!\!/}
\newcommand {\MP} {{\on{MP}}} 

\newcommand {\MU} {{\on{MU}}}
\newcommand {\MUG}[1] {{\MU_G^*(#1)}}

\newcommand {\orb} {{\on{orb}}}
\newcommand {\pf}[1] {{\pi_!^{#1}(1)}}
\newcommand {\PT} {Pontrjagin-Thom }
\newcommand {\pt} {\on{pt}}
\newcommand {\QQ} {\mathbb{Q}}
\newcommand {\res}[2] {\on{res}\arrowvert^{#1}_{#2}}
\newcommand {\RG} {R(G)}
\newcommand {\Sj} {{\Sigma_j}}

\newcommand {\SK} {{{\mathcal S}_{\Kn}}}
\newcommand {\Sl} {{\Sigma_l}}
\newcommand {\Sn} {{\Sigma_n}}
\newcommand {\Snborel} {{\on{E\Sn}\times_\Sn X^n}}

\renewcommand {\SS} {{\mathbb{S}^0}}
\newcommand {\SSn} {{\mathbb{S}^n}}
\newcommand {\SSmn} {{\mathbb{S}^{-n}}}

\newcommand {\Sym} {\on{Sym}}
\newcommand {\tensor}{\otimes}
\newcommand {\T} {T}

\newcommand {\thsG} {{\on{EG}\ltimes_G}}
\newcommand {\TT} {{\T\in\mathcal{T}}}
\newcommand {\TTp} {{\T\in\mathcal{T}_p}}
\newcommand {\td} {{\on{Td}}}
\newcommand {\tdtop} {{\on{Td}_{top}}}
\newcommand {\Tp} {{T_{p^k}}}
\newcommand {\tr}[1] {\on{Trace}(g | #1)}

\newcommand {\Zh} {{\ZZ^h}}
\newcommand {\ZZ} {\mathbb Z}

\newcommand {\zs} {{G/H_+}}
\newcommand {\lps} {[\! [}
\newcommand {\rps} {]\! ]}

\newcommand {\ps}[1] {\lps #1\rps}

\usepackage [all]{xy}
\SelectTips{cm}{12}
\title{Orbifold genera, product formulas and power operations}
\author{Nora Ganter}
\date {\today}
\begin{document}
\begin{abstract}
We generalize the definition of orbifold elliptic genus, and
introduce orbifold genera of chromatic level $h$, using $h$-tuples
rather than pairs of commuting elements. We show that our genera are
in fact orbifold invariants, and we prove integrality results for
them. If the genus arises from an $\hinf$-map into the
Morava-Lubin-Tate theory $E_h$, then we give a formula expressing the 
orbifold genus of the symmetric powers of a stably almost complex
manifold $M$ in terms of the genus of $M$ itself.
Our formula is the $p$-typical
analogue of the Dijkgraaf-Moore-Verlinde-Verlinde formula for the
orbifold elliptic genus \cite{DMVV}.
It depends only on $h$ and not on the genus. 
\end{abstract}
\maketitle
\tableofcontents
\section{Introduction}
This paper aims to provide a systematic understanding and
homotopy theoretic refinement of the theory of orbifold genera and
product formulas as they arise in string theory 
(cf.\ \cite{DMVV}, \cite{Dijkgraaf}).
\subsection{Product formulas}
The most general and famous of these is probably a formula by
Dijkgraaf, Moore, Verlinde and Verlinde expressing the orbifold
elliptic genus of the symmetric powers of an almost complex manifold
$M$ in terms of the elliptic genus of $M$ itself \cite{DMVV}:
\begin{eqnarray}\label{DMVV-formula}
  \sum_{n\geq 0}\phi_{\on{ell},\on{orb}}(\Mn)t^n &=&
  \prod_{\stackrel{m\geq1,n\geq0,}{l\in\ZZ}}\(\frac{1}{1-t^mq^ny^l}\)^{c(mn,l)}\\
  \notag
  &=& \exp\left[\sum_{m\geq1}V_m(\phi_{\on{ell}}(M))t^m\right].  
\end{eqnarray}
Here
$$
  \phi_{\on{ell}}(M) = \sum_{\stackrel{n\geq0,}{l\in\ZZ}}c(n,l)q^ny^l
$$
is the two-variable elliptic genus or ``equivariant $\chi_y$-genus of
the loop space'' of $M$. Its definition can be found in \cite{EOTY},
\cite{Hoehn} or 
\cite{HBJ}. The orbifold version $\phi_{\on{ell},\on{orb}}$ is defined in
\cite{Borisov:Libgober:elliptic}.
The second equation in (\ref{DMVV-formula}) is due to 
Borcherds\footnote{More precisely, Borcherds' computation in the proof
  of the product formula for the $j$-function
  \cite[Lemma 7.1]{Borcherds:monstrous-moonshine} goes through for
  Jacobi forms,  
  if the Hecke operators are replaced by the $V_n$ defined in  
  \cite[I.4.2 (7)]{Eichler:Zagier}.},
and the $V_m$ are a type of Hecke operators acting on $q$-expansions of
Jacobi forms.
Borcherds proved this equality in the context of his proof of the
Moonshine conjectures, and
the right-hand side of (\ref{DMVV-formula}) is often referred to as a 
{\em Borcherds lift} of $\phi_{\on{ell}}(M)$.

We shall show that the $p$-typical analogue of (\ref{DMVV-formula})
arises from a natural equation of cohomology operations in elliptic
cohomology by specializing to the elliptic cohomology of a point.
Thus our work adds to the evidence that elliptic cohomology has a role to
play in the connection between Moonshine and string theory.

Formula (\ref{DMVV-formula}) has been studied by algebraic geometers
\cite{Borisov:Libgober:elliptic} as well as algebraic
topologists \cite{Tamanoi:equivariant-morava-K-theory},
\cite{Tamanoi:banff?}. 
\subsection{The Ando-French definition of orbifold genus}
In \cite{Ando:French}, Ando and
French explain how to fit the notion of orbifold (elliptic) genus into
the framework of equivariant elliptic cohomology. 
The version of equivariant elliptic cohomology they 
choose to work with is Borel equivariant Morava $E$-theory $E_2$.
We explain a slight generalization of Ando and French's definition of
orbifold genus: 
For homotopy theorists, a genus typically is 
a natural transformation from a cobordism theory to another cohomology
theory, applied to a point. If the target of this natural
transformation is a form of elliptic cohomology, for instance $E_2$,
the genus is called an elliptic genus. 
%
%
Let $G$ be a finite group.
For each natural number $h$ there is a Morava $E$-theory
$E_h$, 
and an element $\chi$ of $E_h\BG$ can be viewed as a class function on
$h$-tuples of commuting $p$-power order elements of $G$ (cf.\ Section
\ref{HKR-Sec}). 
We write ${\mathcal N}^{U,G}_*$ for the 
bordism ring of compact, closed, smooth $G$-manifolds
with a complex structure on their stable normal bundle.
We write $\MU^G_*$ for the coefficients of the complex cobordism
spectrum $MU_G$. 
There is a \PT map from $\mathcal N^{U,G}_*$ to $\MU^G_*$ and a
completion map from $\MU^G_*$ to $\MU_*(\BG)$.
For details the reader is referred to Section \ref{equivariant}. 
\begin{Def}
\label{orb-Def}
  Let $\phi$ be a map of ring spectra from $\MU$ 
  into $E_h$. 
  We define the Borel
  equivariant version $\phi_G$ of the genus $\phi$ as the composite
  $$
    \phi_G\negmedspace :{\mathcal N}^{U,G}_*\longto \MU^G_* \longto
    \MU_*(\BG) 
    \stackrel{\phi}{\longrightarrow}
    E_{h*}(\BG) 
  $$
  $$
    (M\acted G)\phantom{XXXXX}\longmapsto\phantom{XXXXX}  \phi_G(M),
  $$ 
  where the first two maps are the \PT map and the completion map.
  We define the orbifold genus associated to $\phi$ to be
  $$
    \phi_{\on{orb}}(M\acted G) := \GG{G}\sum_{\alpha}\(\phi_G(M)\)(\alpha),
  $$
  where the sum runs over all $h$-tuples of commuting elements of
  $p$-power order in $G$. 
\end{Def}
Note that instead of $\MU$ we could have used any of the classical Thom
spectra $\on{MSpin}$, $\on{MO}$, $\on{MSp}$, $\on{MU\langle n\rangle}$,
$\on{MO\langle n\rangle}$ etc.

We shall prove an integrality result about Definition \ref{orb-Def}
and show that it defines in fact an orbifold invariant. 
While most of the literature on orbifold genera is highly
computational, these proofs work on a conceptual level.
They rely on deeply homotopy theoretic
properties of the $\Kn$-local categories, suggesting that stable
homotopy theory provides a good framework for the study of orbifold
phenomena.
\subsection{Power operations}
The left-hand side of the DMVV formula (\ref{DMVV-formula})
involves an object well known
to topologists: the assignment
$$
  M\longmapsto\sum_n(\Mn)t^n
$$
is what is called the {\em total power operation} in cobordism (of a
point). 
The right-hand side of the DMVV formula is a function in
$\phi_{\on{ell}}(M)$ which takes sums into products. Total power operations
also have this property.
Thus it is a natural question to ask whether
formulas like the DMVV formula simply reflect the fact that a
natural transformation $\phi$ preserves power operations.
Such a natural transformation that preserves power operations is 
called an $\hinf$-map.
We shall show that
any $\hinf$-map from a cobordism theory into $E_h$ has a
DMVV-type formula for the induced (orbifold) genus. 
\subsection{The chromatic picture}
  In the case $h=1$, the cohomology theory $E_1$ is $p$-completed
  $K$-theory. The standard example of a genus into (non-$p$-completed)
  $K$-theory is the 
  Todd genus. There is an equivariant version of $K$-theory; here 
  $K_G(\pt)$ is the representation ring $R(G)$, 
  and the character of a representation $\rho$ is the 
  class function
  $$\chi(g)=\on{Trace}(\rho(g)).$$ 
  Definition \ref{orb-Def} can then be formulated without the
  $p$-power order part, and it becomes the definition of the {\em topological}
  Todd genus\footnote{
    In the literature this turns up as the Euler characteristic of the complex
    space $M/G$ \cite{Atiyah:Segal} or the topological Euler characteristic of
    the orbifold $M\mmod G$ \cite{Dijkgraaf}.} 
  of the orbifold $M\mmod G$:
  \begin{Def}\label{top-Todd-Def}
    $$\td_{\on{top}}(M\acted G) := \frac{1}{|G|}\sum_{g\in
    G}\tr{\td_G(M)}.$$ 
  \end{Def}
  In the case $h=2$ Definition \ref{orb-Def} is (up to the factor
  $\GG{G}$) the 
  Definition \cite[6.1]{Ando:French}. Ando and French show that this
  definition is a $p$-typical analogue of the definition of orbifold
  elliptic genus discussed in the literature
  \cite{Borisov:Libgober:elliptic}, \cite{DMVV}.  
  Thus our point of view fits the {\em orbifold elliptic genus} (as defined
  by Ando and French) and its product formula into a common picture
  with the {\em topological Todd genus} of an orbifold and its
  product formula, i.e., the former
  is exactly the chromatic level two analogue of the latter.
\subsection{Statement of results}
  \label{top-Todd-is-integral-Sec}
  A priori $\tdtop$ appears to take values in $\CC$. Note however that
  $$
    \GG G\sum_{g\in G}\tr{-}
  $$
  equals the inner product with the trivial representation. 
  This shows that $\tdtop$ 
  takes integral values. In a similar way, using an inner product
  defined by Strickland, we will prove the following proposition (cf.\
  Corollary 
   \ref{integrality-Cor}): 
\begin{Prop}\label{integrality-introduction-Prop}
  The orbifold genus $\phi_{\on{orb}}$ takes values in
  $E^0_h$.
\end{Prop}
  Definition \ref{orb-Def} is formulated in terms of the $G$-space
  $\MG$ rather than the orbifold $M\mmod G$. It is a
  non-trivial fact that $\phi_{\on{orb}}$ depends only on the
  orbifold, and not on its presentation (cf.\ Theorem
  \ref{forb-is-independent-Thm}): 
\begin{Thm}
  Let $M$ be a compact complex manifold acted upon by a finite group $G$, let
  $N$ be a compact complex manifold acted upon by a finite group $H$,
  and assume 
  that the orbifold quotients $M\mmod G$ and $N\mmod H$ are
  isomorphic as (tangentially almost) complex orbifolds. Then
  $$
    \forb{M\acted G} = \forb{N\acted H}.
  $$
\end{Thm}
The analogues of the Dijkgraaf-Moore-Verlinde-Verlinde formula for our
orbifold genera are given by the following theorem (cf.\ Theorem
\ref{DMVV-Thm}): 
\begin{Thm}\label{main-Thm}
  For any $\hinf$-map $\phi$ from a Thom spectrum into $E_h$ there is a
  formula 
  $$
    \sum_{n\geq 0}\phi_{\on{orb}}(M\mmod\Sn)t^n = \exp\left[\sum_{k\geq
    0}T_{p^k}(\phi(M))t^{p^k}\right]. 
  $$
\end{Thm}
%
There are two side results that are hopefully of independent interest
to homotopy theorists. We 
obtain an explicit formula for the Strickland inner product in
Morava $E$-theory (cf. Corollary \ref{inner-product-Cor}): 
\begin{Prop}
  The Strickland inner product in Morava $E$-theory theory is
  \begin{eqnarray*}
    E^0_h(\BG)\tensor E^0_h(\BG) &   \longto   & E_h^0 \\
    \chi\tensor\xi               & \longmapsto & 
              \GG G\sum_{\alpha}\chi(\alpha)\xi(\alpha),
  \end{eqnarray*}
  where the sum is over all $h$-tuples of commuting elements of
  $p$-power order.
\end{Prop}
The main step in proving that $\phi_{\on{orb}}$ is an orbifold invariant
is a theorem about (equivariant) Spanier-Whitehead duals and the Borel
construction (cf.\ Theorem \ref{Borel-Thm}): 
\begin{Thm}
  Let $G$ be a finite group.
  Then there is an isomorphism of functors from the category of
  finite $G$-CW-spectra to the $\Kn$-local category
  $$
    \Borel D_G(-)\cong D(\Borel (-)),
  $$
  where $D_G$ denotes the $G$-equivariant dual and $D$ denotes the
  dual in the $\Kn$-local category.
\end{Thm}
In \cite{Ando:thesis}, Ando defines Hecke operators on the Morava
$E$-theories, generalizing
the (stable) Adams operations in $K\hat{{}_p}(X)$. 
It was pointed out by Atiyah and Tall \cite{Atiyah:Tall} that the
Adams operations in $K$-theory can be defined using the 
more general theory of $\lambda$-rings due to Grothendieck
\cite{Grothendieck}. 
%
In Section \ref{Atiyah-Tall-Sec},
we offer an Atiyah-Tall-Grothendieck type definition for the $T_{p^k}$,
which generalizes Ando's 
definition to any $\Kn$-local $\hinf$-spectrum. 
\subsection{Acknowledgments}
I want to thank Mike Hopkins, who, after we had a rough start,
turned out to be a wonderful advisor and very
inspiring teacher. Most of
this work was guided by Mike's ideas, and it was a great experience to
work with him. Just as important for my work were the many
mathematical discussions with Haynes Miller. I am deeply indebted to
him for his generosity with his time and for his encouraging support.

I am very grateful to Chris French, Matthew Ando and Charles
Rezk for generously sharing their unpublished work with me and for
long and inspiring conversations. Many thanks also go to Neil
Strickland who never got tired of answering my emails.

I had numerous helpful conversations with Alex Ghitza, Mike Hill,
Johan de Jong and Lars Hesselholt.
Plenty of helpful suggestions from the referee greatly improved the
clarity of this paper.
Andr\'e Henriques, Mike Hill and Alex Ghitza read
earlier drafts of this thesis and offered useful comments.
It is a pleasure to thank them all.

This research is the author's PhD thesis. It was partially supported by
a Walter A.~Rosenblith  
fellowship and by a dissertation stipend from the German Academic 
Exchange Service (DAAD).
\part{The topological Todd genus and its product formula}
\label{top-Todd-Part}
Throughout this paper, the discussion of the
(topological) Todd genus and $K$-theory will
serve as a model for our study of elliptic genera and elliptic cohomology.
\section{The Todd genus from the point of view of stable homotopy
theory}\label{Todd-Sec}
In this section, we recall Conner and Floyd's definition of the Todd 
genus\footnote{Conner and Floyd attribute many of the results mentioned
  here to Atiyah and Hirzebruch \cite{Atiyah:Hirzebruch} and Dold
  \cite{Aarhus}.} 
\cite[I]{Conner:Floyd}. 
\subsection{The Conner-Floyd map}
In \cite{ABS}, 
Atiyah, Bott and Shapiro construct $K$-theory
Thom classes $u_{\on{ABS}}$ for complex vector bundles.
Conner and Floyd \cite[p.29]{Conner:Floyd} show that 
giving $K$-theory Thom classes for complex vector bundles is 
equivalent
to giving a map of spectra
$$\longmap{\td}{\MU}{K}$$
(denoted $\mu_c$ by Conner and Floyd).
On $\MU(n)$, this map is given 
by
$$u_{\on{ABS}}(\gamma_{\on{univ}}^{n})\in [\MU(n),\ZZ\times\on{BU}].$$
The map $\td$ is called the Conner-Floyd map. 
The induced map on homotopy groups,
$$\td_*:=\pi_*(\td ),$$
is the Todd genus. 
\subsection{The push-forward of one}\label{pushforward-Sec}
Assume we are given a multiplicative cohomology theory $E^*(-)$ with
natural Thom classes for complex vector bundles, or 
equivalently, a
map of ring spectra $\map{\phi}{\MU}{E}$. Let $[M]\in\MU_d$, i.e., let
$M$ be a compact closed smooth $d$-dimensional manifold together with
a choice of lift 
$-[\tau]_K\in\widetilde K(M)$ of its stable normal bundle
$-[\tau]\in\widetilde{KO} (M)$. Such manifolds are called ``manifolds with
stably almost complex structure''.
Let $$\longmap{\pi}{M}{\pt}$$ denote the unique map from $M$ to a point.
The following is a slight reformulation of the definition of
``Umkehr'' map along $\pi$ in \cite[pp.40-41]{Dyer}, using the
language of Thom spaces of virtual bundles set up e.g.\ in \cite[IV]{Rudyak}.
\begin{Def} \label{pushforward-Def}
  The push-forward along $\pi$ in $E^*(-)$ is defined by
  $$
    \xymatrix{
      {\pi_!^{\phi}\negmedspace : E^*(M)}\ar[r]^{\cong} &
      {\widetilde E^{*-d}(M^{-\tau})}\ar[r] & 
      {\widetilde E^{*-d}(S^0)\cong \pi_{*+d}(E),} 
    }
  $$
  where the first map is the Thom isomorphism for $-[\tau]_K$ and
  the second map is the Pontrjagin-Thom collapse. 
\end{Def}
\begin{Prop}\label{pushforward-of-1-Prop}
  The genus induced by $\phi$,
  $$
    \longmap{\phi_*}{\MU_*}{E_*}
  $$ 
  sends $[M]\in\MU_d$ to the push-forward of one $\pf{\phi}\in E_d$.
\end{Prop}
\begin{Pf}{}
  The transformation $\phi$ maps Thom classes to Thom classes and thus
  $\pi_!^{\id_\MU}$ to $\pi_!^{\phi}$. Therefore it is sufficient to consider
  the universal case $\phi=\id_\MU$. In this case the statement
  follows directly from the definition of the cobordism Thom classes and
  from the Pontrjagin-Thom construction.  
\end{Pf}
\subsection{The classical definition}
The Riemann-Roch theorem \cite[Dyer]{Aarhus} yields the following
formula for the Todd 
genus\footnote{This expression is the inverse of the one Conner and
  Floyd obtain, 
  because they work with the tangent bundle
  rather than the normal 
  bundle.}:
$$
  \td(M) = \pi^{\td}_!(1) = 
  \int_M\prod_i\frac{1-e^{x_i}}{x_i}.
$$
Here 
the $x_i$ are the Chern roots of the normal bundle $\nu_M$. 
Let $M$ be a compact complex manifold. Then the index theorem implies:
$$
  \td(M) = \sum (-1)^i\dim(H^i(M,\mathscr O_M)),
$$
where $\mathscr O_M$ is the structure sheaf of $M$
(cf.\ \cite[5.4]{HBJ} and \cite[p.542]{Atiyah:Segal}).
\section{The equivariant Todd genus}
\label{equivariant}
All the constructions of the previous section go through
equivariantly.
\subsection{Okonek's equivariant Conner-Floyd maps}
The following proposition and examples are taken from 
\cite[1]{Okonek}\footnote{For an English reference see
  \cite{May:CBMS}. There is a difference between the two: Okonek works
  with tom Dieck's definition of an equivariant cohomology theory
  \cite{tomDieck:Lokalisierung}. In the
  language of \cite{May:CBMS} this is a complex-stable, naive
  $G$-equivariant cohomology theory.  
  }.
\begin{Prop}[Okonek]\label{Okonek-Prop}
  If $E_G^*$ is a multiplicative, $G$-equivariant
  cohomology theory with natural Thom classes for complex $G$-bundles,
  then there is a unique natural, stable transformation
  $$
    \longmap{\phi_G}{\MU_G^*(-)}{E_G^*(-)}
  $$
  of multiplicative $G$-equivariant cohomology theories that takes
  Thom classes to Thom classes.
\end{Prop}
%
Rather than explaining all the concepts in the statement
of the proposition, we state the two examples that are relevant to
us. 
\begin{Exa}\label{Okonek-Borel-Exa}
  For any complex oriented ring spectrum $E$,
  Borel equivariant 
  $E$-cohomology $$E(\borel{-})$$ has natural Thom
  classes for complex $G$-bundles. In this case, $\phi_G$ factors through
  Borel equivariant cobordism. If we let $\phi$ be the orientation of
  $E$, then $\phi$ preserves equivariant Thom classes, so that
  we get 
  $$
    {\phi_G}\negmedspace
    :{\MU_G(-)}\longto\MU(\borel{-})\stackrel{\phi}{\longrightarrow} E(\borel{-}).
  $$
\end{Exa}
\begin{Exa}\label{Okonek-K-Exa}
  (Compactly supported) equivariant $K$-theory has natural Thom
  classes for complex 
  $G$-bundles. We denote the resulting
  equivariant Conner-Floyd maps by
  $$\longmap{\td_G}{\MU_G}{K_G}.$$
\end{Exa}
There is a Pontrjagin-Thom map from the equivariant cobordism
ring ${\mathcal N}^{U,G}_*$ to the coefficient ring $\MU_G^*$, which
in the equivariant
case fails
to be an isomorphism.
\begin{Def} The equivariant Todd genus of an almost complex
  $G$-manifold is defined to be
  $$
    \td_G(M) := \td_{G*}([M]),
  $$
  where $[M]$ denotes the image of $M$ under the $G$-equivariant 
  Pontrjagin-Thom
  map. 
\end{Def}
\subsection{Equivariant push-forward of one}\label{equiv-push-Sec}
  The Thom spectrum $M^\xi$ of a virtual equivariant
  bundle $\xi\in \on{KO}_G(M)$
  and the Thom isomorphism 
  for a choice of stably almost complex
  structure $[\xi]_K\in\widetilde{K}_G(M)$ on $\xi$
  are defined in \cite[X]{Lewis:May:Steinberger}, and 
  Definition \ref{pushforward-Def} goes
  through for an equivariant theory with Thom classes.  
  On the image of the Pontrjagin-Thom map the same argument as in
  Proposition 
  \ref{pushforward-of-1-Prop} shows that
  $$
    \phi_G(M) = \pf{\phi_G} \in E_G(\pt),
  $$ 
  where $\map\pi M{\pt}$ is the unique $G$-map. 

In the case of equivariant $K$-theory, our definition of push-forward is
  equivalent to that of Atiyah and Singer\footnote{More precisely, if
    one replaces 
  $TX$ by $X$ and assumes that all the Thom classes that are needed
  exist, Atiyah and Singer's $\on{ind}_G^X$ becomes our $\pi^X_{!G}$.}
  in \cite{Atiyah:Singer}.
Recall that the correct generalization of the Borel
construction to $G$-spectra is given by the ``twisted half smash
product'' over $G$
$$\thsG-.$$ 
These twisted half smash products were introduced and studied
extensively in \cite{Lewis:May:Steinberger}. 
A summary of their basic properties can be found in \cite[I.1]{BMMS}.
For the suspension
spectrum of a pointed $G$-space $X$, they specialize to the Borel
construction 
$$
  \thsG(\Sigma^\infty X)\cong\Sigma^\infty(\Borel X).
$$ 
In the case of the Thom spectrum $M^{-\tau}$ we have
(cf. \cite[X.6.3]{Lewis:May:Steinberger})
$$
  \thsG (M^{-\tau}) = (\borel M)^{-\borel\tau}.
$$
%
\subsection{An explicit formula for $\on{Td}_G$}
The equivariant Todd genus takes values in the representation ring
$K_G(\pt)=\RG$, and it is classical that a representation 
$V\acted G$ is determined by its character
$$
  g\longmapsto\tr V.
$$
  Using the Riemann-Roch theorem and a Lefschetz fixed point
  formula, 
  Atiyah, Segal and Singer (cf.\ \cite[(2.11)]{Atiyah:Segal} and
   \cite{Atiyah:Singer:III}) prove the following:
  \begin{equation}\label{Atiyah-Segal-Singer-Eqn}
  \tr{\td_{G}(M)} = \int_{M^g}\frac{1}{\prod_\zeta 1-\zeta 
    e^{x_j(N^g_\zeta)}},
   \end{equation} 
  where $M^g$ stands for the $g$-fixed points of $M$, 
  and $N^g$ denotes the normal bundle of $M^g$ in $M$;
  the product runs over all eigenvalues of the action of $g$ on $N^g$
  and over the Chern roots $x_j(N^g_\zeta)$ of the eigenbundles  
  $N^g_\zeta$.

At this point we would like to point out how Definition \ref{orb-Def} 
relates to that of Borisov and Libgober \cite{Borisov:Libgober:elliptic}, 
which looks like the right-hand side of
(\ref{Atiyah-Segal-Singer-Eqn}). Recall that Definition 
\ref{orb-Def} follows the one given by Ando and French, who
generalize the left-hand side of (\ref{Atiyah-Segal-Singer-Eqn}).
Character theory is available in
the context of Ando and French's work, 
but the Riemann-Roch formula is not. However, they
explain in detail how to modify the character theoretic
discussion in \cite{Atiyah:Segal} to bring their definition into a
form that is modulo a Riemann-Roch theorem very similar to Borisov and
Libgober's. Their discussion goes through without changes for our Definition
\ref{orb-Def}. 
\section{Power operations}\label{power-op-Sec}
\subsection{Power operations and $\hinf$-ring spectra}
Let 
$\{E_G\mid G \text{ finite}\}$ 
be a compatible family of equivariant cohomology theories in the sense of
\cite[II.8.5]{Lewis:May:Steinberger}, and write $\EG X$ for $E_G^0(X)$. 
``Compatible'' implies
in particular that for a map $\alpha\negmedspace:H\to G$ and a $G$-space
$X$, we have a restriction map
$$
  \longmap{\res{}\alpha}{E_G(X)}{E_H(X)},
$$
and if $\alpha$ is the inclusion of a subgroup and $X$ a $G$-space we
also have an induction map
$$
  \longmap{\ind HG}{E_H(X)}{E_G(X)},
$$
such that the axioms of a Mackey structure on $E_G$ spelled out in
\cite{tomDieck:Mackey} are satisfied. 
  The author could not find the reference for this fact, so here is a
  short explanation: 
  A compatible family satisfies 
  \begin{equation}\label{compatible-Eqn}
    E_G(G\ltimes_H X) \cong E_H(X)
  \end{equation}
  for $H\sub G$ and any (pointed) $H$-space $X$. 
  If $X$ is already a $G$-space, one has 
  $$
    G\ltimes_H X\cong\zs\wedge X,
  $$
  and the map
  $$
    \longmap{(p_{G/H})_+}{\zs}{\SS}
  $$
  sending all of $G/H$ to the non-basepoint induces $\res GH$, while
  its $G$-equivariant Spanier-Whitehead dual 
  $$
    \longmap{D_G(p_{G/H})_+}{\SS}{\zs}
  $$
  induces $\ind HG$. For arbitrary $\alpha$, the compatibility
  condition does not provide us with an isomorphism
  (\ref{compatible-Eqn}), but with a map from the left to the
  right. Thus if  
  we replace
  $(p_{G/H})_+$ by the co-unit of the adjunction
  $(G\ltimes_\alpha-,\on{forget})$, we can still define $\res{}\alpha$.
  The Mackey criteria
  follow from \cite[XIX.3]{May:CBMS}.
We also ask that our family has unitary, commutative and associative
{\em external products}  
$$
  \longmap{\wedge}{E_G(X)\tensor E_H(Y)}{E_{G\times H}(X\wedge Y)}, 
$$
that are natural in (stable) maps of $X$ and $Y$. Note that
this implies that $\wedge$ also commutes with induction
and restriction maps.
By {\em unitary} we mean that for each $G$, there is an element
$1\in\EG\SS$ with
$1\wedge x = \res{}{\on{pr}_2}x$, where $\on{pr}_2$ is the projection
onto the second factor of $G\times H$. We further ask that $\res{G}{H}1=1$.
\begin{Exa}
  For any $E$, Borel equivariant $E$-cohomology $E(EG\times_G-)$ is an
  example \cite[XXI.1.9]{May:CBMS}. Here the induction maps equal the
  transfer maps
  $$
    \longmap{T_H^G}{\Sigma^{\infty}_+(\borel
    X)}{{\Sigma^\infty_+(\on{EH}\times_HX)}}.
  $$ 
\end{Exa}
\begin{Exa}
  Equivariant $K$-theory is an example, with the induction maps 
  the induced representation.
\end{Exa}
Before we recall the definition of an $\hinf$-structure on $\{E_G\}$,
we need to introduce some notation.
  Let $X$ be a pointed $G$-space. We write
  $$
    (X\acted G)^n\acted\Sn
  \text{\phantom{X}    or \phantom{X}   }
    X^n\acted (G\wr\Sn)
  $$
  for the space $X^{\wedge n}$ acted on by
  $$
    \GS = G^n\rtimes\Sn
  $$
  as follows: $G$ acts on each factor individually, while
  $\Sn$ permutes the factors.
  By abuse of notation, we also write $E_\Sn(X^n)$ for
  $E_{G\wr\Sn}(X^n)$ and in particular $E(X)$ for $E_G(X)$, unless we
  want to emphasize the equivariant situation.
The following definition is essentially \cite[VIII.1.1]{BMMS}.
\begin{Def}\label{hinf-Def}
  An $\hinf$-structure on $E$ is given by a collection of natural maps
  $$
    \longmap{P_n}{E_G(X)}{E_\GS(X^n)}
  $$ 
  called {\em power operations} satisfying the following conditions: 
  \begin{enumerate}
    \renewcommand{\labelenumi}{(\alph{enumi})}
    \item $P_1 = \id$ and $P_0(x) = 1$,
    \item the (external) product of two power operations is 
      $$P_j(x)\wedge P_k(x) = 
        \res{\Sigma_{j+k}}{\Sigma_j\times\Sigma_k} (P_{j+k}(x)),
      $$
    \item the composition of two power operations is
      $$
        P_j(P_k(x)) = \res{\Sigma_{jk}}{\Sigma_k\wr\Sigma_j}\(P_{jk}(x)\),  
      $$
    \item and the $P_j$'s preserve (external) products:
    $$P_j(x\wedge y) = \res{\Sj\times\Sj}\Sj(P_j(x)\wedge P_j(y)),$$
  \end{enumerate} 
  where the restriction is along the map
  $$
    \left[((X\acted G)^2)^j\acted\Sj\right] \longrightarrow \left[(X\acted
    G)^{2j}\acted(\Sj\times\Sj)\right]\cong
    \left[((X\acted G)^j\acted\Sj)^2\right]. 
  $$
\end{Def}
\begin{Rem}\label{Einf-Rem}
Traditionally people formulated this definition only for Borel equivariant
theories. In that case it is a refinement of the notion of ring
spectrum up to homotopy, but 
it is weaker than the notion of
$E_\infty$ or $A_\infty$ structure. 
In the same way our definition is weaker than Greenlees and May's
notion of {\em global $\mathcal I_*$ functor with smash product
  spectrum} in \cite{Greenlees:May:completion}.
More precisely, using the Yoneda lemma one can reformulate Definition
\ref{hinf-Def} 
in terms of maps of
$\GS$-spectra
$$
  \longmap{\xi_n}{E_G^n}{E_\GS}.
$$
The maps $\res GH$ are induced by maps of $G$-spectra
$$
  G\ltimes_H E_H \longto E_G,
$$
and the external product $\wedge$ becomes a map
of $G\times H$-spectra.  
Thus conditions (1)-(4) of the definition translate into homotopy
commutative diagrams of spectra.
A global $\mathcal I_*$ functor with smash product spectrum has such
$\xi_n$, and in that case, the diagrams commute strictly. 
\end{Rem}
%
%
%
\subsection{Total power operations}\label{total-power-Sec}
Let $E$ be an $\hinf$-ring spectrum. It is often convenient to
consider all power operations at once, i.e. the {\em total power
  operation}
$$\longmap{P}{E(X)}{\cOplus E_\Sn(X^n)t^n}$$
which is $P_n$ into each summand. Here we are following the notation
of \cite{Segal:Aspen}: The symbol $\hat{\bigoplus}$ stands for the
infinite product, and
the variable $t$ is a dummy
variable, 
introduced in order to keep track of the ``summand'' and also to
avoid convergence issues later on.
  Note that $$\cOplus E_\Sn(X^n)t^n$$ is a graded ring by 
  $$E_\Sn(X^n)\otimes E_{\Sigma_m}(X^m)\stackrel{\wedge}{\longrightarrow}
  E_{\Sn\times\Sigma_m}(X^{n+m}) 
  \stackrel{\on{ind}}{\longrightarrow}
  E_{\Sigma_{n+m}}(X^{n+m}),$$
  where 
  $$
    \on{ind} = \ind{\Sn\times\Sigma_m}{\Sigma_{m+n}}
  $$      
  (compare \cite{Segal:Aspen}). 
\begin{Prop}[compare {\cite[VIII.1.1]{BMMS}}]\label{properties-of-P_js-Prop}
  We have
  \begin{enumerate}
    \renewcommand{\labelenumi}{(\alph{enumi})}
    \item
    the restriction  of $P_j(x)$ to $E(X^j)$ is 
    $$
      \res\Sj1 P_j(x) = x^{\wedge j},
    $$
    \item the operation $P_j$ applied to $1\in
    E^0(\SS)$ is 
    $$
      P_j(1) = 1_\Sj := 1 \in E_\Sj(\pt),
    $$
    \item\label{P-is-exp} the total power operation takes sums
    into products: 
    $$
      P(x+y) = P(x)\cdot P(y).
    $$
  \end{enumerate}
\end{Prop}
\begin{Pf}{}
  The first two properties are immediate from the definition. The
  proof of property (c) in 
  \cite[II.1.6]{BMMS} and \cite[VII.1.10]{Lewis:May:Steinberger} takes
  place on the level of equivariant spectra.
\end{Pf}
\subsection{Power operations in $K$-theory and cobordism}
In \cite{Atiyah:poweroperations} Atiyah defines power operations for
$K$-theory. In the case of an 
(equivariant) vector bundle $V$ over (a $G$-space) $X$, they are
given by the (external) tensor product 
$$P_n([V]) = [V^{\otimes n}]\in K_\Sn(X^n).$$
%

In \cite{tomDieck:poweroperations} tom Dieck defines power operations
for Borel equivariant cobordism and shows that the Conner-Floyd map is
an $\hinf$-map in the classical (i.e.\ Borel equivariant) sense. 
We prefer to work on the level of
equivariant cobordism $\MUG{-}$ (cf.\ \cite{tomDieck:Bordism}).
In that case, Greenlees and May show that ${\MU_G}$ is a 
``global $\mathcal I_*$ functor with smash product
spectrum''
\cite[5.8]{Greenlees:May:completion}. Thus, 
by Remark \ref{Einf-Rem} it has power operations.
On coefficients\footnote{More precisely: on non-equivariant coefficients
  or on the image of the Pontrjagin-Thom map.}
these power operations in equivariant cobordism are
\begin{eqnarray*}
  P_n\negmedspace :\MU^*(\pt) &   \longto   & \MU^{n*}_\Sn(\pt) \\
  {[M]}       & \longmapsto & [\Mn]. 
\end{eqnarray*}
\begin{Prop}[{compare \cite[(A4)]{tomDieck:poweroperations}}]
  \label{MUG-universal-Prop}
  The $P_n$ are multiplicative with respect to $\wedge$ and 
  compatible with Thom classes in the following sense:
  \begin{equation}\label{A4-Eqn}
    P_n(u_\MU(\xi)) = u_{\MU,\Sn}(\xi^{\oplus n}),
  \end{equation}
  and $\{ MU_G\}$ equipped with the $P_n$ is universal with respect to
  this
  property.  In other words, for any equivariant cohomology theory
  with multiplicative Thom classes 
  for complex $G$-bundles and power operations satisfying
  (\ref{A4-Eqn}) the maps $\phi_G$ of Proposition \ref{Okonek-Prop}
  preserve power operations.
\end{Prop}
\begin{Pf}{}
  If $\xi$ is the universal complex $G$-bundle, Equation
  (\ref{A4-Eqn}) is immediate from the construction; for
  other $\xi$ it follows by naturality. Let now $\{E_G\}$ be an
  equivariant family as in the proposition. Then the proof of
  Proposition \ref{Okonek-Prop} 
  shows that
  the $\phi_G$ preserve power operations.  
\end{Pf}
\begin{Cor}[compare \cite{tomDieck:poweroperations}]\label{C-F-hinf-Cor}
  The map
  $$
    \MUG{-}\longrightarrow\MU^*(\borel{-})
  $$ 
  of Example \ref{Okonek-Borel-Exa}
  and
%
  the equivariant Conner-Floyd-Okonek map 
  $$\longmap{\on{Td}_G}{\MU_G}{K_G}$$
  of Example
  \ref{Okonek-K-Exa} are $\hinf$-maps. 
\end{Cor}
\begin{Pf}{}
  Both maps are defined as examples of the map $\phi_G$,
  which is an $\hinf$-map by 
  the proof of Proposition \ref{MUG-universal-Prop}.
\end{Pf}
  Let $E$ be an $\hinf$-spectrum with compatible Thom classes as in 
  Proposition \ref{MUG-universal-Prop}, and let $V$ be a complex
  $d$-dimensional $G$-representation. Then $E_G$ 
  comes equipped with natural isomorphisms
  $$E_G^0({\mathbb S}^{2d}\wedge
  X)\stackrel{\cong}{\longrightarrow}E^{0}(V^{c}\wedge 
  X),$$
  where $V^c$ denotes the one point compactification of $V$ (cf.\
  \cite[2.1]{Greenlees:May:completion}). 
  This becomes important when
  we want to extend our power operations to $$E^{-2d}(X) = E^0({\mathbb
  S}^{2d}\wedge X),$$ because $({\mathbb S}^{2d})^n\acted\Sn$ is an
  equivariant sphere. In the situation of the proposition we can
  follow \cite{Greenlees:May:completion} to
  extend the power operations to 
  $$\longmap{P_n}{E^{2d}(X)}{E^{2nd}(X^n)}.$$
\subsection{Internal power operations}\label{internal-power-operations-Sec}
\label{Atiyah-Sec}
We can always compose the power operation $P_n$ with the pullback
along the diagonal map of $X^n$
$$\longmap{\Delta_n^*}{E_\Sn(X^n)}{E_\Sn(X)}.$$
Since the action of $\Sn$ on $X$ is trivial, the target of this map
often turns out to be 
$${E_\Sn(\pt)\otimes E(X)}.$$
We might want to compose further with a map
$$E_\Sn(\pt)\longto E^0$$
in order to obtain 
operations\footnote{In the literature (e.g.\ \cite{Ando:thesis}) these
  compositions are often referred to as power operations and $P_n$ is
  then called ``total power operation''. We follow the convention to
  call them {\em internal power operations}, since they actually act
  on $E(X)$.}
acting on $E(X)$.
In the case of $K$-theory, $E_\Sn(\pt)$ is the representation ring
$R(\Sn)$. 
The following two examples from
\cite{Atiyah:poweroperations} are important to us.
\begin{Exa}
  Atiyah's definition of the Adams operations is
  $$\psi_n(x) = \on{Trace}(c_n|\Delta_n^*P_n(x)),$$
  where $c_n$ is a cycle of length $n$.
\end{Exa}
\begin{Exa}\label{symmetric-powers-Exa}
  The operations $\sigma_n$ are defined by
  \begin{eqnarray*}
    \sigma_n(x)   & := &
    \frac{1}{n!}\sum_{g\in\Sn}\tr{\Delta_n^*P_n(x)}\\  
    & = & \langle \Delta_n^*P_n(x),1\rangle_\Sn.
  \end{eqnarray*}
  If $x=[V]$ is the class of a vector
  bundle $V$, then 
  $$\sigma_n(x) = [\Sym^n(V)]$$
  is represented by the $n^{th}$ symmetric power of $V$, since in this
  case the inner product with $1$  
  counts the multiplicity of the trivial representation as a
  summand of
  $$[V^{\otimes n}\acted\Sn] = \Delta_n^*P_n(x)\in R(\Sn)\otimes K(X).$$ 
\end{Exa}
\begin{Def}\label{symmetric-powers-Def}
  We write
  $$S_t(x) := \sum_{n\geq 0}\sigma_n(x)t^n$$
  for the total symmetric power. 
  In other words, $S_t$ is the composite
  $$S_t\negmedspace : K(X)\stackrel{P}{\longto}\cOK\longto
  \cOplus R(\Sn)\otimes K(X)t^n\longto K(X)\ps{t},$$
  where on the $n^{th}$
  summand the second map is pullback along the diagonal 
  and the third map is the inner product with $1_\Sn$. 
\end{Def} 
In \cite{Ando:thesis}, Ando generalizes
Atiyah's work to cohomology theories with \HKR
  character theory, as we will recall in Section \ref{Matthew-Sec}. 
\section{Multiplicative formulas for the Todd genus}\label{Todd-product-Sec}
The following is a reformulation of the second statement of
Corollary \ref{C-F-hinf-Cor}:
%
\begin{Cor}\label{Todd-H-infty}
  The square
  $$\xymatrix{
    {\MU(X)}\ar^{\td}[0,3]\ar_{\text{\raisebox{2.7ex}{$P_{\MU}$}}}[d]& & &
    K(X)\ar^{\text{\raisebox{1.7ex}{$P_K$}}}[d] \\ 
    \cOMU \ar^{\hat\Oplus\td_\Sn}[0,3]  & & &\cOK
  }$$
  commutes.
\end{Cor}
It follows immediately that the equivariant Todd genera $\td_\Sn(M^n)$
are determined by the Todd genus of $M$, and moreover that
the expression is exponential in $\td(M)$. More precisely,
specializing to the case where $X$ is a point results in the
following corollary:
\begin{Cor}\label{toddformula}
  Let $M$ be an almost complex manifold. Then we have the following
  equation in the ring $\cOplus R(\Sn)t^n$:
  $$\sum_{n\geq 0} \td_\Sn (M^n)t^n = \(\sum_{n\geq 0}1_{\Sn}t^n
  \)^{\td(M)},$$
  where $1_\Sn\in R(\Sn)$ denotes the trivial representation.
\end{Cor}
\begin{Pf}{}
  By Corollary \ref{Todd-H-infty} we have 
  $$\sum_{n\geq 0}\td_\Sn(M^n) t^n = P_K(\td(M)).$$
  By Proposition \ref{properties-of-P_js-Prop} (c), $P_K$
  takes sums into products.
  Since 
  $$\td(M)\in K(\pt) = \ZZ,$$
  this implies
  $$P_K(\td(M)) = P_K(1)^{\td(M)}.$$ 
  Now
  $P_n(1)$ is the trivial representation of $\Sn$
  (compare Proposition \ref{properties-of-P_js-Prop} 
  (b)). Therefore,
  $$P_K(1)=\sum_{n\geq 0}1_\Sn t^n.$$ 
\end{Pf}
%
%
As a further consequence of Corollary \ref{Todd-H-infty}, we obtain the
multiplicative formula for the topological Todd genus
\cite{Dijkgraaf}:
\begin{Cor}\label{toptoddformula-Cor}
We have
$$
  \sum_{n\geq 0}\on{Td}_{\on{top}}(M^n\acted\Sn)t^n = 
  \(\frac{1}{1-t}\)^{\td(M)} =
  \exp \left[ \sum_{n\geq 1} \frac{\psi_n(\on{Td}(M))}{n}t^n \right].
$$
\end{Cor}
This is the chromatic level one analogue of the DMVV formula
(\ref{DMVV-formula}).
%
\medskip

\begin{Pf}{}
  We have   
  \begin{eqnarray}
  \label{Todd-symmetric-Eqn}
    \sum_{n\geq 0}\tdtop(\Mn)t^n & = & 
    \sum_{n\geq 0}\frac{1}{n!}\sum_{g\in\Sn}\tr{\td_\Sn(M^n)}t^n \\ 
\notag 
   &=&\sum_{n\geq 0}\frac{1}{n!} \sum_{g\in\Sn}\tr{P_n(\td(M))}t^n \\
\notag
    &=& S_t(\td(M)),
  \end{eqnarray}
  where the first equation is the definition of $\td_{\on{top}}$, the second
  equation is Corollary \ref{Todd-H-infty}, and the third equation is
  Definition \ref{symmetric-powers-Def} with $X$ the one point space. 

The first identity of the corollary now follows exactly like Corollary
  \ref{toddformula} from the 
fact that $S_t$ is exponential.
We thank Charles Rezk for reminding us of the well-known equation
$$
  S_t(x) = \exp \left[ \sum_{n\geq 1} \frac{\psi_n(x)}{n}t^n \right].
$$
Together with (\ref{Todd-symmetric-Eqn}) this proves 
$$
  \sum_{n\geq 0}\on{Td}_{\on{top}}(M^n\acted\Sn)t^n = 
  \exp \left[ \sum_{n\geq 1} \frac{\psi_n(\on{Td}(M))}{n}t^n \right].
$$
\end{Pf}

\part{The orbifold elliptic genus and other higher chromatic relatives
  of $\tdtop$} 
\label{higher-chromatic-Part}
The methods of Part \ref{top-Todd-Part} appear to be specific to
equivariant $K$-theory: We use the inner product of two
representations, symmetric powers of vector bundles, and
evaluation of characters at group
elements. 
%
Our discussion in the higher chromatic case relies on the fact
that character theory as well as inner products have been defined in
much greater generality.
Firstly, for $E$ a suitable $\Kn$-local cohomology theory, e.g.\ 
Morava $E$-theory $E_h$, and $\chi$ an element of $\EBG$,
Hopkins-Kuhn-Ravenel theory defines evaluation of $\chi$ at $h$-tuples
of commuting $p$-power order elements of $G$ (cf.\ Section \ref{HKR-Sec}).
Secondly, Strickland has defined inner products
$$
  b_G\negmedspace :\EBG\otimes\EBG\longto E^0
$$
in any $\Kn$-local cohomology theory $E$
%
(cf.\ Section \ref{inner-product-Sec}).
If \HKR theory applies and $E^0$ is torsion free
they satisfy the formula
\begin{equation}
  \label{inner-product-Eqn}
  b_G(\chi,\xi) = \GG G\sum_{\alpha}\chi(\alpha)\xi(\alpha),  
\end{equation}
where the sum runs over all $h$-tuples of commuting $p$-power order
elements of $G$
(cf.\ Corollary \ref{inner-product-Cor}).

In Sections \ref{Matthew-Sec} and \ref{orb-genera-Sec}, we recall how to
use \HKR character theory to define
orbifold genera $\phi_\orb$ and Hecke operators in Morava
$E$-theory. We also generalize the definition of symmetric powers to
operations in Morava $E$-theory.
If $\phi$ is an $\hinf$-map
we prove a DMVV-type product
formula for $\phi_\orb$ (cf.\ Section \ref{proof-of-main-thm-Sec}). 
The formula (\ref{inner-product-Eqn}) implies integrality statements
for $\phi_\orb$ and the symmetric powers (cf.\ Corollaries
\ref{integrality-Cor} and \ref{sym-integrality-Cor}). 
Another consequence of (\ref{inner-product-Eqn}) is the key fact that
the map 
$$
  \chi\longmapsto\GG G\sum_{\alpha}\chi(\alpha)
$$
is induced by a map in the $\Kn$-local category. 
It will play a central role in Section \ref{McKay-Sec}, where we prove
that $\phi_{\on{orb}}(\Mn)$ 
does not depend on the representation of the orbifold $M\mmod G$. 
It will also allow us to generalize the definitions of symmetric powers
and Hecke operators to any $\Kn$-local $\hinf$-spectrum $E$ (cf.\
Sections \ref{generalized-symmetric-powers-Sec} and \ref{Atiyah-Tall-Sec}).
\section{\HKR theory}
This section recalls some results from \cite{HKR}.
The reader can find a nice and short
introduction to \HKR character theory in
\cite[5]{Ando:French}, 
also see \cite[8]{Rezk:logarithms}.
\subsection{Even periodic ring spectra and formal groups}
We keep our paper in the language of even periodic ring spectra,
because all our examples are of this kind. 
This section is a short reminder of their definition and
properties. For details see \cite{AHS}.
\begin{Def}
  An even periodic ring spectrum is a spectrum $E$ such that the
  graded coefficient ring $E^*$ is concentrated in even degrees and
  $E^2$ contains a unit. 
\end{Def}
No choice of this unit is specified.
  In the context of even periodic ring spectra it is often convenient
  to replace the complex cobordism spectrum $\MU$ by its two-periodic
  version 
  $$
    \MP := \bigvee_{j\in\ZZ}\Sigma^{2j}\MU.
  $$
  Note that $\MP$ is the Thom spectrum of $\ZZ\times\on{BU}$.
  For even periodic $E$
  the Atiyah-Hirzebruch spectral sequence for $E^*(\CC P^n)$
  collapses, and the system 
  $$  
    E^*(\CC P^n)\longleftarrow E^*(\CC P^{n+1})
  $$
  is Mittag-Leffler, such that $E^*(\CC P^\infty)$ becomes
  non-canonically isomorphic to $E^*\ps x$. As usual\footnote{Cf.\
    \cite{Adams}, \cite{Rudyak}.}
  a good choice of such an $x$ gives rise to $E$-theory Chern classes, and
  to a formal group law $F$ over $E_*$ describing the first Chern
  class of the tensor product of line bundles
  $$
    c_1(L_1\tensor L_2) = c_1(L_1)+_Fc_1(L_2).
  $$
  The advantage of working with even periodic $E$ is that rather than
  speaking about formal group laws one can use the language of formal
  groups: For such $E$ the map 
  $$\CC P^\infty\times\CC P^\infty\longto\CC P^\infty$$ classifying the
  tensor product of line bundles makes the formal spectrum
  $\on{spf}E^0(\CC P^\infty)$ into an (affine one-dimensional) formal
  group scheme, and choosing a coordinate for this formal group is
  equivalent to specifying a map of ring spectra
  $$
    \on{MP}\longto E.
  $$
  We do not make much use of these concepts, but we use several
  results whose proofs rely on a deep understanding of the way these
  formal groups come into the picture. For the moment it
  is enough to remember that an even periodic ring spectrum $E$ has
  somehow a formal group attached to it.
%
\subsection{Morava $E$-theories}\label{HKR-Sec}
We now explain which spectra we can work with.
\begin{Def}\label{HKR-Def}
  Let $E$ be an even periodic ring spectrum with associated
  formal group $F$. We say that {\em $E$ has a \HKR theory} if 
  \begin{enumerate}
  \renewcommand{\labelenumi}{(\alph{enumi})}
  \item $E^0$ is local with maximal ideal $\mathfrak m$, and complete
    in the $\mathfrak m$-adic topology,
  \item the graded residue field $E^0/\mathfrak m$ has characteristic
    $p>0$, 
  \item $p^{-1}E^0$ is not zero,
  \item the mod $\mathfrak m$ reduction of $F$ has height $h<\infty$
    over $E^0/\mathfrak m$.
  \end{enumerate}
\end{Def}
%
Hopkins, Kuhn and Ravenel give a list of examples satisfying the
conditions of this definition. One of these examples is in addition an
$\hinf$-spectrum and the interplay between \HKR theory and the
$\hinf$-structure is well understood. This is the reason why it
becomes our favorite example: 
\begin{Exa}[Lubin-Tate cohomology/Morava $E$-theory]
  Consider the graded ring
  $$
    E_* := \mathbb{WF}_{p^h}\ps{u_1,\dots,u_{h-1}}[u^{\pm 1}],
  $$
  where $u_i$ has degree zero, $u$ has degree $2$, and $\mathbb{W}k$
  denotes the ring of Witt vectors of the field $k$.
  There is a cohomology theory called {\em Lubin-Tate cohomology} or
  {\em Morava
  $E$-theory}, which has $E_*$ as coefficients.   
  On a finite complex $X$, it is given by 
  $$
    E_h^*(X) = \MU^*(X)\tensor_{\MU^*}E^*, 
  $$
  where the map $\MU_*\to E_*$ classifies the universal deformation of
  the Honda formal group law.
  The construction of this cohomology theory goes back many years,
  a published account can be
  found in\footnote{Rezk 
    omits a subtlety in his exposition: He proves that
    $E_*$ is Landweber exact over $\on{BP}$, obtaining a homology
    theory. Via Spanier-Whitehead duality this becomes a cohomology
    theory on finite complexes as described. 
    The phantom discussion in \cite{Hovey:Strickland} proves
    that it is (uniquely) represented by a ring spectrum.} 
  \cite{Rezk:HM}.
\end{Exa}
\subsection{$h$-tuples of commuting elements}
Just as classical characters of $G$ are class functions on $G$, \HKR
characters are class functions on $h$-tuples of commuting $p$-power
order elements of $G$. This section is a short reminder of the basic
definitions concerning such $h$-tuples; we will give a more detailed
discussion of the case $G=\Sn$ in Section \ref{h-tup-Sec}.   
Since $G$ is finite, the set of all \htpps of $G$ can be identified with
$$\Hom(\ZZ_p^h,G).$$
The group $G$ acts on this set by conjugation:
$$g\gg g^{-1} = (gg_1g^{-1},\dots,gg_hg^{-1}).$$
\begin{Def}
  Let $\alpha$ be an $h$-tuple of commuting elements (of $p$-power
  order) of $G$. The conjugacy class $[\alpha]_G$ of $\alpha$ is
  defined to be the orbit of $\alpha$ in
  $\Hom(\ZZ^h,G)$ 
  (or $\Hom(\ZZ_p^h,G)$ respectively) under this $G$ action.
  The centralizer of $\alpha$ is defined as the stabilizer
  $$
    C_\alpha = C_G(\alpha) := \on{Stab}_G(\alpha) \sub G.
  $$ 
\end{Def}
\begin{Def}
  A function on $\Hom(\ZZ_p^h,G)$ is called a class function if it is
  invariant under conjugation by elements of $G$.
\end{Def}
\subsection{\HKR characters}
Let $E$ be a spectrum with \HKR theory, 
let $G$ be a finite group, 
let $\chi$ be an element of $E^0(BG)$ and
let $\alpha$ be an $h$-tuple of commuting $p$-power order elements of
$G$, where $h$ is as in Definition \ref{HKR-Def}. Then 
Hopkins, Kuhn and Ravenel define a ring $D$ and an evaluation map
\begin{eqnarray*}
  {\on{eval_\alpha}}\negmedspace : \EBG & \longto & {D} \\
    \chi & \longmapsto & \chi(\alpha).
\end{eqnarray*}
For our purposes it is not important what the ring $D$ is or how
$\on{eval_\alpha}$ is defined, but for completeness, we
recall their definitions:
Let $D_n$ be the ring
$$
  D_n := E^0(B(\ZZ/p^n\ZZ)^2)/\text{(annihilators of nontrivial Euler classes)},
$$
then 
$$
  D=\Colim_{n}D_n
$$
is the colimit over the maps induced by
$$
  \ZZ/p^{n+1}\ZZ\longto\ZZ/{p^n}\ZZ.
$$
Since $G$ is finite, any $\alpha\in\Hom(\ZZ_p^h,G)$ factors through some
$\alpha_n\in\Hom((\ZZ/p^n\ZZ)^h,G)$, and 
$$\eval(\chi):=\alpha_n^*(\chi)\in D$$ is independent of the choice of $n$.
  
We will use the fact that $D$ is independent of the group $G$ and that a
fixed $\chi\in \EBG$ defines a class function on 
the set
of $h$-tuples of commuting $p$-power order elements of $G$. This is
the sense in which $\chi$ is a character.
The maps $\eval{}$ are analogues of the $\tr{-}$ maps in
representation theory.
%
%
The following is a corollary of \cite[Thm C]{HKR}.
\begin{Thm}\label{C-Thm}
Let $E$ be a ring spectrum with \HKR theory.
      An element $\chi$ of $\frac{1}{p}E^0(BG)$ is uniquely
      determined by the class function it defines. 
\end{Thm}
We also need the \HKR analogue of the formula for the 
character of an induced representation
\cite[p.30]{Serre:representations}. 
\begin{Thm}[{\cite[Thm D]{HKR}}]
  \label{D-Thm}
  Let $H\sub G$ be a subgroup, and let $\alpha$ be an $h$-tuple of
  commuting $p$-power order elements in $G$. We have
  \begin{eqnarray*}
    (\ind HG(\chi))(\alpha) & = & \GG{H}\sum_{\stackrel{g\in G\mid g\alpha
    g^{-1}}{\text{maps to H}}} \chi(g\alpha g^{-1}). 
  \end{eqnarray*}
\end{Thm}
\subsection{Ando's generalization of Atiyah's work}\label{Matthew-Sec}
The original reference for this Section is \cite{Ando:thesis}, see also
\cite{Ando:Duke}. 
Let $E$ be an $\hinf$-ring spectrum with \HKR theory. 
For simplicity we assume a K\"unneth isomorphism for the symmetric groups,
i.e., we ask that $E^0(\on{B\Sn})$ be free
of finite rank over $E^0$ and that $E^1(\on{B\Sn})=0$. 
\begin{Exa}({\cite[3.3]{Strickland:symmetric},\cite{Ando:Duke}})
  The Morava $E$-theories $E_h$ satisfy all the above conditions.
\end{Exa}
For such spectra, Ando defines internal power operations.
The examples relevant to us are the analogues of the examples 
in Section \ref{Atiyah-Sec}. 
\begin{Def}
  Let $\alpha$ be an \htpp of $\Sn$. Define
  $$
    \psi_\alpha \negmedspace : E(X) \longto D\tensor E(X)  
  $$ 
  as the composition 
  $$
    E(X)\stackrel{P_n}{\longrightarrow} E(\Snborel)
    \stackrel{\Delta_n^*}{\longrightarrow} E(\on{B\Sn}\times X)
    \stackrel{\cong}{\longleftarrow} E(\on{B\Sn})\otimes E(X)
    {\longrightarrow} D\tensor E(X),
  $$
  where $\Delta_n$ denotes the diagonal map of $X^n$, and the last arrow
  sends $\chi\tensor x$ to $\chi(\alpha)\tensor x$.
\end{Def}
  Let $\alpha$ be as above. Then $\alpha$ makes $\{1,\dots,n\}$ into a
  $\ZZ^h_p$-set. Conversely, a finite $\ZZ^h_p$-set $A$ determines an
  $h$-tuple of $p$-power elements $\alpha$ of some symmetric group up
  to conjugacy (cf.\ 
  Section \ref{h-tup-Sec}). We sometimes write $\psi_A$ for
  $\psi_\alpha$. 
\begin{Def}\label{Hecke-Def}
  The Hecke operators in Morava-Lubin-Tate theory are defined as
  $$
    T_{p^k}(x) := \frac{1}{p^k}\sum_{\stackrel{\TTp}{|\T|=p^k}}\psi_\T(x),
  $$
  where the sum is over all isomorphism classes of
  transitive $\ZZ^h_p$-sets of order
  $p^k$. 
  It is proved in \cite{Ando:thesis} that
  these $T_{p^k}$ are additive operations
  $$
    \longmap{\Tp}{E_h(X)}{E_h(X)}.
  $$   
\end{Def}
  Note that on $E_1=K\hat{{}_p}$, these Hecke operators are the stable Adams
  operations: 
  $$
    T_{p^k} = \frac{\psi_{p^k}}{p^k}.
  $$
\begin{Def}\label{higher-symmetric-powers-Def}
  Let $E^0$ be torsion free. We define the analogues of the symmetric 
  powers as
  $$
    \sigma_n(x) := 
    \frac{1}{n!}\sum_{\alpha}\psi_\alpha(x),
  $$
  where $x\in E(X)$, and
  this time the sum runs over all $h$-tuples $\alpha$ of
  commuting elements of $p$-power order in $\Sn$.
  We write $S_t$ for the ``total symmetric power'' as above.
\end{Def}
  It is immediate from \cite[5.5]{Ando:French} that the operation
  $\sigma_n$ takes values in $\GG\Sn E^0(X)$, but
  it is a non-trivial fact that it takes values in
  $E^0(X)$. We postpone the proof to Section
  \ref{integrality-Sec}. 
  As in the case of $K$-theory, $S_t$ turns out to take sums into
  products. We will give a more general definition of the $\sigma_n$ and
  prove this exponential property in Section 
  \ref{generalized-symmetric-powers-Sec}. 
\section{Generalized orbifold genera}\label{orb-genera-Sec}
Recall from Definition \ref{orb-Def} that if $M$ is a stably almost
complex oriented $G$-manifold, and $$\longmap\phi\MU{E_h}$$ is a
complex orientation of Morava $E$-theory, then the orbifold genus
$\forb\MG$ is defined by the formula
$$
  \forb\MG = \GG{G}\sum_{\alpha}\(\phi_G(M)\)(\alpha),
$$
where $\phi_G$ is the Borel equivariant version of $\phi$, and the sum
runs over all $h$-tuples of commuting elements of 
$p$-power order in $G$.     
In this section, we generalize the definitions of $\phi_\orb$ and
of $\sigma_n$, using maps in the $\Kn$-local category.
As a corollary of these new definitions,
we obtain the promised integrality statements.
%
\subsection{The $\Kn$-local categories}\label{Kn-Sec}
Let $H_*(-)$ be a generalized homology theory. Recall from
\cite{Bousfield}
that there is a category $\mathcal S_H$, called the
$H$-local (stable homotopy) category, and a functor
$$
  \longmap\gamma{\mathcal S}{\mathcal S_H},
$$
which is left-universal with respect to the property that it takes
$H_*$-isomorphisms into isomorphisms.
When it is clear that we are working in $\mathcal S_H$, we
will often omit $\gamma$ from the notation.
Like the stable homotopy category $\mathcal S$ itself, $\mathcal S_H$ is
a triangulated category with 
a compatible closed symmetric monoidal structure. In other words,
it has a symmetric monoidal structure $-\wedge-$ with unit
$S=\gamma(\SS)$ and function objects (``internal $\hom$'s'') $F(-,-)$,
such that
$$
  \Hom(X\wedge Y,Z) = \Hom(X,F(Y,Z)), 
$$
and these data are compatible with the triangulated structure in an
appropriate 
sense\footnote{The details can be found in
  \cite[A.2]{Hovey:Palmieri:Strickland}.}.
The localization functor $\gamma$ preserves the triangulated structure
as well as the monoidal structure and
its unit, but does not in general preserve function objects\footnote{Cf.\
  \cite[3.5.1]{Hovey:Palmieri:Strickland}.}.  
There is one important class of function objects preserved by
$\gamma$, which is going to play a role for us:
Write 
$$
  DX := F(X,S)
$$
for the dual of $X$.
\begin{Thm}[{\cite[III.1.6]{Lewis:May:Steinberger}}]
  \label{strong-dual-Thm}
  Let $X$ and $Y$ be objects of a closed
  symmetric monoidal category, and assume that
  there are maps
  $$
    \longmap{\alpha}{S}{X\wedge Y} \quad\text{ and   }\quad
    \longmap{\beta}{Y\wedge X}{S} 
  $$
  such that the composites
  $$(\id\wedge\beta)\circ(\alpha\wedge\id) \negmedspace :
    X\cong S \wedge X \longrightarrow X\wedge Y\wedge X
    \longrightarrow X\wedge S\cong X
  $$
  and
  $$(\beta\wedge\id)\circ(\id\wedge\alpha) \negmedspace :
    Y\cong Y\wedge S\longrightarrow Y\wedge X\wedge Y
    \longrightarrow S\wedge Y\cong Y
  $$
  are the respective identity maps. Then the adjoint
  $\map{\beta^\sharp}{Y}{DX}$ is an isomorphism.
\end{Thm}
An object $X$ for which such $Y$, $\beta$ and $\alpha$ exist is called
{\em strongly dualizable}. It comes with an isomorphism $X\to DDX$.
Since $\gamma$ preserves the monoidal
structure, Theorem \ref{strong-dual-Thm} implies
that $\gamma$ also preserves strong dualizability and
strong duals.
\begin{Def}
  Let $E$ be a spectrum such that any map that becomes an
  isomorphism under $H_*(-)$ also becomes 
  an isomorphism under $E^*(-)$. Then $E$ is called an $H$-local
  spectrum.
\end{Def}
If $E$ is $H$-local, $E^*(-)$ is a well-defined functor on the category
$\mathcal{S}_H$. 
The following theorem seems to be well-known to homotopy 
theorists\footnote{To the author's knowledge there is no published 
  account of it. In the case that $E$ is Morava-Lubin-Tate cohomology
  it is proved in  
  \cite[5.2]{Hovey:Strickland}, for Noetherian $E^0$ a written account
  is available from \cite{Strickland:communication}, in the generality
  it is stated here I learned it from Michael Hopkins.}:
\begin{Thm}\label{HKR-Kn-Thm}
  Let $E$ be a cohomology theory with level $h$ \HKR
character theory. Then $E$ is local with respect to the Morava
$K$-theory $K(h)$. 
\end{Thm}
These {\em Morava $K$-theory} homology theories $\Kn_*(-)$
were first constructed by Baas and Sullivan and first used by Morava.
Today their definition can be found in \cite{Rudyak} or \cite{EKMM}.
  The functor $\gamma$ has a fully faithful right-adjoint $J$,
  whose image is the (full) subcategory of $\Kn$-local spectra, 
  and it
  is customary to think of $\SK$ as embedded into
  $\mathcal S$ via $J$. 
  This point of view is not helpful for our
  purposes, and we stick to the language of localized categories. The
  difference is mainly in notation: 
  Write $L_\Kn$ for the composite $J\circ\gamma$. The functor $J$ does
  not preserve the monoidal structure. Thus,  
  where we write
  $$\gamma(X)\wedge\gamma(Y)\quad \text{ or }\quad X\wedge Y$$
  for the smash product in $\SK$, others write
  $$L_\Kn(L_\Kn X\wedge L_\Kn Y),$$
  and similarly we write $\SS$ or $\gamma(\SS)$ for $L_\Kn\SS$.
\subsection{Strickland inner products}\label{inner-product-Sec}
This section recalls some of the concepts and results in
\cite{Strickland:inner:product}.
Let $\mathcal C$ be an additive closed symmetric monoidal
category. We use the notation of the previous section and write $\tau$
for the twist map $X\wedge Y\to Y\wedge X$.
We fix the assumption on $\mathcal C$ that
every object is strongly dualizable.
\begin{Def}\label{Frobenius-Def}
  A {\em Frobenius object} in $\mathcal C$ is an object $A$ equipped
  with maps 
  $$S\stackrel{\eta\phantom{i}}{\longto} A,\quad A\wedge
  A\stackrel{\mu\phantom{i}}{\longto} 
  A,\quad A\stackrel{\eps\phantom{i}}{\longto} S,\quad\text{and }
  A\stackrel{\psi\phantom{i}}{\longto} 
  A\wedge A$$ 
  such that
  \begin{enumerate}
    \renewcommand{\labelenumi}{(\alph{enumi})}
    \item $(A,\eta,\mu)$ is a commutative and associative monoid,
    \item $(A,\eps,\psi)$ is a commutative and associative
      co-monoid,
    \item we have $\psi\circ\mu = (1\wedge\mu)\circ(\psi\wedge 1)$.
  \end{enumerate}
\end{Def}
\begin{Lem}[{\cite[3.9]{Strickland:inner:product}}]\label{inner-product-Lem}
  If $(A,\eta,\mu,\psi,\eps)$ is a Frobenius object in $\mathcal
  C$ then $b:=\eps\mu$ defines an inner product on $A$ in the
  following sense:  
  \begin{enumerate}
    \renewcommand{\labelenumi}{(\alph{enumi})}
    \item $b$ is symmetric, i.e., $b\circ \tau=b$, and
    \item $b$ is non-degenerate, i.e., the adjoint
      $\map{b^\sharp}{A}{DA}$ is an isomorphism.
  \end{enumerate}
\end{Lem}
Let $G$ be a finite group or groupoid, and let $\BG$ denote its Borel
construction.
We write $\BG_+$ for the $\Kn$-local suspension
spectrum of the Borel contruction of $G$:
$$\gamma(\Sigma^\infty_+\BG).$$ 
%
Let $\map\delta G{G\times G}$ denote the inclusion of the diagonal,
and write $\psi$ for $B\delta_+$ and $\mu$ for the transfer map 
$$
  \longmap{\mu=T\delta}{\BG_+\wedge\BG_+}{\BG_+}.
$$
Let $p_G$ be the unique map from $G$ to the trivial group, and write
$\eps$ for $\on{Bp_G}_+$. Let $\beta$ be the composite
$\eps\circ\mu$. 
In the following, $D$ will denote the dual in the $\Kn$-local category.
\begin{Thm}[{\cite[8.7,3.11,8.2,8.5]{Strickland:inner:product}}]
  \label{Neil-Thm}  
  In the $\Kn$-local category, $\beta$ is an inner product on
  $\BG_+$. Let
  $$
    \eta\negmedspace :\SS=D\SS\longto D\BG_+\stackrel\cong\longleftarrow\BG_+
  $$ 
  be the composite of $D\eps$ with $(\beta^\sharp)^{-1}$.
  Then $(\BG_+,\mu,\eta,\psi,\eps)$ is a Frobenius object in the
  $\Kn$-local category.
\end{Thm}
From now on, let $E$ be an even periodic 
$K(h)$-local spectrum. Then unreduced $E$-cohomology of the space
$\BG$ is the same as $E$-cohomology of the spectrum $\BG_+$, and we
write $E^0(\BG)$ for both.
Let $m$ be the composite
$$ 
  m\negmedspace : E^0(\BG)\otimes E^0(\BG)\longto E^0(\BG\times\BG) \longto
  E^0(\BG)
$$
of the K\"unneth map with $\psi^*$. Then $(m,\eps^*)$ is the standard
ring structure on $E^0(\BG)$, and $$b_G:=\eta^*\circ m$$ defines a
symmetric bilinear form on $E^0(\BG)$.
Assume that $E^0(\BG)$ has finite rank over $E^0$. In this case, the
K\"unneth map becomes an isomorphism over $\QQ$, and the map 
$$
  \mu^*=\ind{\delta}{}
$$
defines a comultiplication on $\QQ\tensor E^0(\BG)$,
$$
  \longmap{\mu^*}{E^0(\BG)}{E^0(\BG\times\BG)\cong_\QQ
  E^0(\BG)\tensor E^0(\BG)}.
$$
\begin{Cor}\label{Neil-Cor}
  If $E^0(\BG)$ has finite rank over $E^0$,
  the maps $m$, $\eps^*$, $\mu^*$ and $\eta^*$ make $\QQ\tensor
  E^0(\BG)$ into a 
  Frobenius object in the category of $\QQ\tensor E^0$-modules.
\end{Cor}
%
%
%
Note that $b_G$ is defined integrally, but that there it might not satisfy the
  non-degeneracy condition (b) of Lemma \ref{inner-product-Lem}.
  Note also that the augmentation map $\eta^*$ is the same as the inner
  product with $1$:  
  $$
    b_G(\chi,1) = \eta^*(\chi\cdot 1) =  \eta^*(\chi).
  $$
The proof of Frobenius reciprocity
\cite[p.25]{Strickland:inner:product} goes through (integrally) in our
situation:
\begin{Prop}
  Let $\map iHG$ be an inclusion of finite groups. Then we have
  $$ 
    b_G(\ind HG\chi,\xi) = b_H(\chi,\res GH\xi).
  $$
\end{Prop}
\begin{Pf}{}
  Let $\xi=1$. We have 
  \begin{eqnarray*}
    \eta^*_G(\ind HG\chi) & = & \eta^*_G(Ti)^*(\chi)\\
     & = &
     \eta^*_G\circ(\beta_G^\sharp)^*
     \circ(D\on{Bi}_+)^*\circ((\beta_H^\sharp)^{-1})^*(\chi)\\   
     & = & (D\eps_G)^*\circ ((\beta_G^\sharp)^{-1})^*\circ
(\beta_G^\sharp)^*
     \circ(D\on{Bi}_+)^*\circ((\beta_H^\sharp)^{-1})^*(\chi)\\
     & = & (D\eps_H)^*\circ ((\beta_H^\sharp)^{-1})^*(\chi)\\
                       & = &  \eta^*_H(\chi),
  \end{eqnarray*}
where the first equation is the definition of $\ind HG$, the second
equation is \cite[8.5]{Strickland:inner:product}, the third and the
last equation follow from 
the definition of $\eta$, and the fourth equation follows from 
$p_G\circ i=p_H$ and the definition of $\eps$.
  Let
  now $\xi$ be arbitrary. Let $\map{j}{H}{H\times G}$ denote
  the diagonal inclusion. Note that
  $$
    \res{}j =  \res{}{\delta_H}\circ(\id\times\res GH).
  $$
  The proof of \cite[8.5]{Strickland:inner:product} implies
  $$
    \ind HG\circ\res{}j = \res{}{\delta_G}\circ(\ind HG\times\id).
  $$
  Combining these three equations, we obtain
  \begin{eqnarray*}
    b_H(\chi,\res GH\xi) & = &
              (\eta^*_H\circ\res{}{\delta_H})(\chi,\res GH\xi)
              \\ 
    & = & \eta_G^*\circ\ind HG\circ\res{}{\delta_H}(\chi,\res GH\xi) \\
    &=& \eta^*_G\circ
         \ind HG\circ\res{}{j}(\chi,\xi)\\
    & = &  \eta^*_G\circ\res{}{\delta_G}(\ind HG\chi,\xi)\\
    & = & b_G(\ind HG\chi,\xi),
  \end{eqnarray*}
where the first and the last equation are the definitions of $b_H$ and
$b_G$. 
\end{Pf}
\subsection{Integrality theorem}\label{integrality-Sec}
The goal of this section is to prove the following proposition:
\begin{Prop}\label{augmentation-Prop}
  Over $\QQ$ the augmentation map $\map{\eta^*}{E^0(BG)}{E^0}$ is 
  $$
    (\eta^*\tensor\QQ)(\chi) = \GG G\sum_\alpha\chi(\alpha).
  $$
\end{Prop}
\begin{Cor}
  If $E^0$ is torsion free, the right hand side defines a map 
  $$E^0(\BG)\longto E^0.$$
\end{Cor}
%
%
%
\begin{Cor}\label{integrality-Cor}
  The orbifold genus 
  $$\phi_{\on{orb}}(\MG) = \GG G\sum_\alpha(\phi_G(M))(\alpha)$$ 
  of Definition \ref{orb-Def} 
  takes values in $E_h^0$.
\end{Cor}
\begin{Cor}\label{sym-integrality-Cor}
  The symmetric powers 
  $$
    \sigma_n(x) := 
    \frac{1}{n!}\sum_{\alpha}(\Delta_n^*P_n(x))(\alpha)
  $$
  of Definition
  \ref{higher-symmetric-powers-Def} take values in $E^0(X)$. 
\end{Cor}
\begin{Pf}{of Proposition \ref{augmentation-Prop}} 
  Let $a$ denote the map
  $$
    a\negmedspace : \chi\longmapsto\GG G\sum_\alpha\chi(\alpha).
  $$
  We need to show that $a$ 
  is a unit of $\QQ\tensor\mu^*$. Since units of co-multiplications are
  uniquely determined, this implies that $a$ is equal to $\eta^*\tensor\QQ$.
  We first compute $\mu^*=\ind{\delta}{}$ in terms of Hopkins-Kuhn-Ravenel
  characters. Let 
  $$
    (\alpha,\beta) := ((a_1,b_1),\dots,(a_h,b_h))
  $$
  be an $h$-tuple of commuting elements of
  $p$-power order in $G\times G$.
  Then by Theorem \ref{D-Thm}
  $$
    \(\ind{\delta}{}(\chi)\)
    (\alpha,\beta) = \GG G 
    \sum_{\stackrel{(s,t)}{s^{-1}\alpha s = t^{-1}\beta t}} 
    \chi(s^{-1}\alpha s).
  $$
  Thus, counting the pairs $(s,t)$ and taking into account that 
  $\chi(s^{-1}\alpha s)=\chi(\alpha)$, we have
  $$
    \(\mu^*(\chi)\)(\alpha,\beta) = \GG G \sum_{s\in G}
    \sum_{\stackrel{t\in G}{s^{-1}\alpha s =
    t^{-1}\beta t}} \chi(\alpha) =
    \begin{cases}
       |C_\alpha|\cdot \chi(\alpha) & \alpha\sim_G\beta \\
       0                       & \text{else}. 
    \end{cases}
  $$
  We are now ready to prove that $a$ is a unit of $\mu^*\tensor\QQ$,
  i.e.\ that the equality 
  $$
    (\id_{\LT{\BG}}\tensor a)\circ\mu^* = \id_{\LT{\BG}} 
  $$
  holds over $\QQ$. 
  By Theorem \ref{C-Thm} it suffices to show that both sides define
  the same class function.
  Write 
  $$
    \xi := (\id\tensor a)\circ\mu^* (\chi).
  $$
  We have
  \begin{eqnarray*}
    \xi(\alpha) & = & \GG
    G\sum_\beta(\mu^*(\chi))(\alpha,\beta) \\
    & = & \GG G\sum_{\beta\in[\alpha]_G}|C_\alpha|\cdot\chi(\alpha) \\
    & = & \chi(\alpha).
  \end{eqnarray*}
\end{Pf}
As a further corollary of Proposition \ref{augmentation-Prop} we
obtain the formula for the Strickland inner product mentioned in the
introduction: 
\begin{Cor}\label{inner-product-Cor}
  Let $E$ be a cohomology theory with Hopkins-Kuhn-Ravenel theory,
  and assume that $E^0$ is torsion free. 
  Then the Strickland inner product on $E^0(\BG)$ is described by the
  formula
  $$b_G(\chi,\xi) = \GG G\sum_\alpha \chi(\alpha)\xi(\alpha).$$
\end{Cor}
\begin{Pf}{}
  We have
  $$
    b_G(\chi,\xi) = (\eta^*\circ m)(\chi,\xi) = \eta^*(\chi\cdot\xi)
    =\GG G\sum_\alpha \chi(\alpha)\xi(\alpha),  
  $$  
  where the first equation is the definition of $b_G$, the second is
  the fact that $m$ is the standard multiplication on $E^0(\BG)$,
  and the last equation follows from Proposition
  \ref{augmentation-Prop}, since $E^0$ is torsion free.
\end{Pf}
\subsection{Generalized orbifold genera}\label{general-def-Sec}
We are now ready to give our most general definition of orbifold genus. 
Recall that Definition \ref{orb-Def} requires an even periodic 
cohomology theory $E$ with level $h$ \HKR theory, and that any such $E$ is 
$\Kn$-local. 
Proposition \ref{augmentation-Prop} motivates the following definition:
\begin{Def}\label{generalized-genus-Def}
  Let $E$ be an even periodic $K(h)$-local ring spectrum, and 
  let $\map\phi\MU E$ be a map of ring
  spectra. Let $G$ be a finite group, and let $\phi_G$ be the Borel
  equivariant genus associated to $\phi$ as in Definition \ref{orb-Def}.
  We define the orbifold genus $\phi_\orb$ of stably almost complex
  $G$ manifolds as the composition 
  \begin{eqnarray*}
    {\phi_{\on{orb}} := \eta^*\circ\phi_G}:{{{\mathcal N}^{U,G}_*}}&\longto&{E_*},
  \end{eqnarray*}
  where $\eta$ is the map of Theorem \ref{Neil-Thm}.
\end{Def}
Instead of $\MU$ we could have used any of the classical Thom
spectra $\on{MSpin}$, $\on{MO}$, $\on{MSp}$, $\MU\langle n\rangle$,
$\on{MO}\langle n\rangle$ etc.
  In the case $E=E_h$,
  Definition \ref{generalized-genus-Def} specializes  by Proposition
  \ref{augmentation-Prop} to the Ando-French 
  Definition \ref{orb-Def}.
\subsection{Generalized symmetric powers}
\label{generalized-symmetric-powers-Sec} 
Recall the definitions of symmetric powers in $K$-theory  
(Example \ref{symmetric-powers-Exa}) and in $E$-cohomology, where $E$ is 
an 
$\hinf$-spectrum with \HKR theory and a K\"unneth isomorphism for the 
symmetric groups
(Definition \ref{higher-symmetric-powers-Def}).
Proposition \ref{augmentation-Prop} motivates the
following generalization of these definitions:
\begin{Def}\label{generalized-symmetric-powers-Def} 
  Let $E$ be an even periodic $\Kn$-local $\hinf$-ring spectrum.
  Let $X$ be a space with basepoint or a spectrum.
  We define 
  the $n^{th}$ symmetric power in $E(X)$ by
  $$\sigma_n := (\eta_{\Sn}\wedge\id_X)^*\circ\Delta_n^*\circ P_n$$
  and the total symmetric power by
  $$
    S_t\negmedspace : E(X) \longto \cOplus E(\on{E\Sn}_+\wedge_\Sn X^n) t^n
    \longto \cOplus E(\on{B\Sn}_+\wedge X) t^n
    \longto E(X)\ps{t} 
  $$
  $$S_t(x) = \sum_{n=0}^\infty \sigma_nt^n(x),$$
where the first map is the total power operation, and on the $n^{th}$ summand, 
the second map is pullback along the diagonal of $X^n$, while the third map
is pullback along $\eta_\Sn\wedge\id_X$.
\end{Def}
In the situation of Definition \ref{higher-symmetric-powers-Def}, 
Theorem \ref{HKR-Kn-Thm} implies that $E$ is $\Kn$-local, and the
two definitions agree by Proposition \ref{augmentation-Prop}.
Note that Definition \ref{generalized-symmetric-powers-Def} does not
require a K\"unneth condition like the one in Definition
\ref{higher-symmetric-powers-Def}. 
We are now going to show that $S_t$ is exponential.
Recall from Section \ref{total-power-Sec} 
that $$\cOplus E^0(\on{E\Sn}_+\wedge_\Sn X^n)t^n$$
is a ring, where multiplication is defined using the transfer maps
$\ind{\Sigma_{m}\times\Sn}{\Sigma_{n+m}}$,
$$
  E_\Sn(X^n)\tensor E_{\Sigma_m}(X^m)\longto
  E_{\Sn\times\Sigma_m}(X^n\times X^m)\longto E_{\Sigma_{n+m}}(X^{n+m}).
$$
\begin{Lem}\label{evaluationismultiplicative}
The map
$$
  \(\sum_{n\geq 1} (\eta_\Sn\wedge\id_X)^*\)\circ\(\sum_{n\geq
  1}\Delta_n^*\) \negmedspace : 
  \cOplus E^0(\on{E\Sn}_+\wedge_\Sn X^n)t^n\longto E^0(X) \ps{t}
$$
is a map of rings.
\end{Lem}
\begin{Cor}
  The total symmetric power $S_t$ takes sums into products.
\end{Cor}
\begin{Pf}{}
  Total power operations take sums into products (cf.\ Propostition
  \ref{properties-of-P_js-Prop} (c)), and
  $S_t$ is defined as a total power operation followed by the ring map
  of the lemma.
\end{Pf}
\begin{Pf}{{of Lemma \ref{evaluationismultiplicative}}}
  Note first that the target of $\sum_{n\geq 0}\Delta_n^*$,
  $$\cOplus E^0(\on{B\Sn}_+\wedge_\Sn X)t^n,$$
  also carries a ring structure: 
  the multiplication is defined by 
  $$
    E_\Sn(X)\tensor E_{\Sigma_m}(X)\longto
    E_{\Sn\times\Sigma_m}(X\times X)
    \stackrel{\Delta_2^*}{\longrightarrow}
    E_{\Sn\times\Sigma_m}(X) \longto E_{\Sigma_{n+m}}(X),
  $$  
  where $E_G$ denotes Borel equivariant $E$-cohomology and
  the last map is again $\ind{\Sn\times\Sigma_m}{\Sigma_{n+m}}$.
  We have
  $$
    \Delta_{n+m} = (\Delta_{n}\times\Delta_m)\circ\Delta_2
  $$
  as maps of $\Sigma_{n+m}$-spaces, and $\ind HG$ is natural in maps
  of $G$-spaces. Therefore the map
  $$
    \sum_{n\geq 1}\Delta_n^*
  $$
  is a ring map. It remains to show that
  $$
    \sum_{n\geq 1} (\eta_\Sn\wedge\id_X)^*
  $$
  is a map of rings.
  Recall that $$\longmap{\eps_G}{\BG_+}{\SS}$$ is
  $B(-)_+$ applied to the unique map from $G$ to the trivial
  group, and that $\eta_G = D\eps_G$. Thus 
  $$
    \eps_{G\times H}=\eps_G\wedge\eps_H \quad\text{ and }\quad
    \eta_{G\times H}=\eta_G\wedge\eta_H.
  $$
  Together with Frobenius reciprocity, this implies
  \begin{eqnarray*}
    \eta^*_{\Sigma_{n+m}}\circ\ind{\Sn\times\Sigma_{m}}{\Sigma_{n+m}}
    & = &
    (T_{\Sn\times\Sigma_{m}}^{\Sigma_{n+m}}\circ\eta_{\Sigma_{n+m}})^* \\
    & = & \eta^*_{\Sigma_n\times\Sigma_m}\\
    & = & (\eta_{\Sigma_n}\wedge\eta_{\Sigma_m})^*.
  \end{eqnarray*}
  This proves the lemma if $X$ is a point.
  Together,
  \begin{eqnarray*}
    (\eta_{\Sigma_{n+m}}\wedge\id_X)^*\circ
    (T_{\Sn\times\Sigma_{m}}^{\Sigma_{n+m}}\wedge\Delta_2)^*  
    & = & ((\eta_\Sn\wedge\eta_{\Sigma_{m}})\wedge\Delta_2)^*\\
    & = &
    (\id_\SS\wedge\Delta_2)^*\circ(\eta_\Sn\wedge\eta_{\Sigma_{m}}\wedge\id_{X\wedge
    X})^*. 
  \end{eqnarray*}
\end{Pf}
\section{The orbifold genus $\phi_{\on{orb}}$ as orbifold
  invariant}\label{McKay-Sec} 
Let $G$ and $H$ be finite groups acting on smooth manifolds $M$ and
$N$ respectively. In this section we recall the notion of tangentially
almost complex structure and prove the following theorem:
\begin{Thm}\label{forb-is-independent-Thm}
  If the orbifold quotients $M\mmod G$ and $N\mmod H$ are isomorphic
  as orbifolds with almost complex structure, then 
  $$\forb{M\acted G}=\forb{N\acted H}.$$
\end{Thm}
Note however that our definition of $\phi_{\on{orb}}$ only makes sense
for orbifolds which can be represented as a global quotient
$M\mmod G$ 
by a finite group $G$. 
\begin{Rem}
  For Borisov and Libgober's definition of the orbifold elliptic genus
  the analogous statement is a consequence of the McKay correspondence,
  proved in \cite{Borisov:Libgober:McKay}. 
\end{Rem}
We use the following facts about orbifolds:
an isomorphism of orbifolds
$$
  M\mmod G\cong N\mmod H
$$
induces an isomorphism of (real or complex) equivariant $K$-groups
$$
  K_G(M)\cong K_H(N)
$$
and a homotopy equivalence of Borel constructions 
$$
  \on{EG}\times_GM\simeq \on{EH}\times_H N,
$$
such that the following diagram commutes
$$
  \xymatrix{{K_G(M)}\ar[r]\ar@{=}[d]&{K(\on{EG}\times_GM)}\ar@{=}[d]\\
  {K_H(N)}\ar[r]&{K(\on{EH}\times_HN).}}
$$
Here the horizontal arrows are the completion maps, and we will use the
notation 
$$
  \longmap{\on{Borel}}{\Korb{M\mmod G}}{K(\on{Borel}(M\mmod G))}
$$
if we want to emphasize its independence of the representation of the
orbifold. For more background on orbifolds, we refer the reader to 
\cite{Moerdijk}.
\subsection{Tangentially almost complex structures}\label{almost-cx-Sec}
Recall that a (stably) almost complex structure on a $G$-manifold $M$ is a
choice of lift $-[\tau]_K\in\widetilde K_G(M)$ of the stable normal bundle
$-[\tau]\in\widetilde{KO}_G(M)$. 
The tangent vector bundle is a well-defined orbifold notion
\cite{Satake}, but 
$$
  \widetilde K_G(M) = \on{coker}(K_G(\pt)\longto K_G(M)) 
$$
is not, since there is no fixed group $G$.
We can however define
$$
  \text{\v K}_{\on{orb}}{\underline X} := \on{coker}(K_{\orb}(\pt)\longto 
  \Korb{\underline X}), 
$$
for an arbitrary orbifold $\underline X$, and similarly 
\v K$O_{\on{orb}}$. 
\begin{Def}[{compare \cite[XXVIII.3.1]{May:CBMS}}]
  A tangentially almost complex structure on an orbifold $\underline
  X$ is a choice of lift $[\tau]_{\text{\v K}}\in\text{\v
  K}_{\on{orb}}(\underline X)$ of $[\tau]\in
  \text{\v K}O_{\on{orb}}(\underline X)$.
\end{Def}
%
%
%
%
The reduced completion map
$$
  \text{\v K}_{\on{orb}}(M\mmod G)\longto\tilde K(\on{Borel}(M\mmod G))
$$
sends a tangentially almost complex structure on $M\mmod G$ to a lift
$$-[\on{Borel}(\tau)]_K\in\tilde K(\on{Borel} M\mmod G)$$
in such a way that the Borel-equivariant Thom isomorphism
$$
  E^0(\borel M)\longrightarrow E^{-d}(\thsG M^{-\tau})
$$
defined by $-[\tau]_{\check{K}}$ agrees with the non-equivariant Thom
isomorphism 
$$
  E^0(\on{Borel}(M\mmod G))\longrightarrow E^{-d}((\on{Borel}(M\mmod
  G))^{-\on{Borel}(\tau)}) 
$$
defined by $-[\on{Borel}(\tau)]_K$
(compare Section \ref{equiv-push-Sec}).
In particular, this Thom isomorphism is independent of the
representation of the orbifold.
\subsection{The genus $\phi_\orb$ as orbifold invariant}
Recall from Section \ref{general-def-Sec} 
that $$\forb \MG\in E_d$$ is 
the image of one under the composite
\begin{equation}
  \label{Th-PT-eps-Eqn}
  E^0(\borel M)\longrightarrow 
  E^{-d}(\thsG M^{-\tau})\stackrel{PT}{\longrightarrow}
  E^{-d}\BG\stackrel{\eta^*}{\longrightarrow}E^{-d}\SS,
\end{equation}
where the first map is the Thom isomorphism defined by
$-[{\tau}]_{\check{K}}$, 
the second map is the 
Pontrjagin-Thom collapse and the third map is pullback along the
map $\eta$ of Section \ref{inner-product-Sec}.
We have just seen that the Thom isomorphism is independent of the
representation of $M\mmod G$, as long as we fix a tangentially almost
complex structure.
The second and third map in (\ref{Th-PT-eps-Eqn}) clearly depend on
the representation of the orbifold, because $\BG$ does. 
However, we will show that the composition
\begin{equation}
\label{eta-PT-Eqn}
  \on{Borel}(PT)\circ\eta :\SS\longrightarrow\thsG M^{-\tau}
\end{equation}
in the $\Kn$-local category
is independent of the representation.
More precisely, we will prove the following theorem:
\begin{Thm}\label{dual-independent-Thm}
  The map (\ref{eta-PT-Eqn}) is the Spanier-Whitehead dual in the $\Kn$-local 
  category of the map
  $$
    (\borel M)_+ \longrightarrow\SS,
  $$
  sending $\borel M$ to the non-basepoint of $\SS$.
\end{Thm}
As a corollary, we obtain Theorem \ref{forb-is-independent-Thm}.
\medskip

\begin{Pf}{{of Theorem \ref{forb-is-independent-Thm}}}
  The map of Theorem \ref{dual-independent-Thm}
  is independent of the representation of $M\mmod G$.
  Therefore, so are
  the maps in (\ref{eta-PT-Eqn}) and (\ref{Th-PT-eps-Eqn}).
\end{Pf}

The remainder of this section is devoted to the proof of Theorem
\ref{dual-independent-Thm}.
\subsection{Borel construction and duality}
%
  We retain the notation of Section \ref{Kn-Sec}
  and write 
  $$
    D_G(-) := F_{G}(-,\SS)
  $$
  for the $G$-equivariant dual (cf.\
  \cite[XVI.7]{May:CBMS}) and 
  $$ 
    D(-) := F_{\mathcal S_\Kn}(-,\SS)
  $$
  for the dual in the $\Kn$-local category.
%
%
\begin{Thm}\label{Borel-Thm}
  Let $G$ be a finite group and let $Y$ be a finite $G$-CW spectrum. 
  There is an isomorphism in the $\Kn$-local category
  \begin{equation}
    \label{Borel-Eqn}
    \thsG(D_G(Y))\longrightarrow D(\thsG Y),
  \end{equation}
  which is natural in $Y$.
\end{Thm}
In order to construct the map in (\ref{Borel-Eqn}), we construct its
adjoint 
\begin{equation}
  \label{adjoint-Eqn}
  \left(\thsG(D_G(Y))\right)\wedge(\thsG Y)\longrightarrow\SS.
\end{equation}
We construct (\ref{adjoint-Eqn}) as a map in the (non-localized)
stable homotopy category, but the adjunction that yields
(\ref{Borel-Eqn}) is in the $\Kn$-local category.
%
  Recall that for a finite group $G$ and $G$-spectra $X$ and $Y$ 
  there is an isomorphism of functors
  $$
    (\thsG X)\wedge(\thsG Y) \cong (\on{EG}\times\on{EG})\ltimes_{G\times
    G}(X\wedge Y).
  $$ 
  The map (\ref{adjoint-Eqn}) is defined as the composite of several
  maps: the first is transfer along the diagonal $\delta$ of $G$ 
  $$ 
    {T\delta}\negmedspace :
    {(\on{EG}\times\on{EG})\ltimes_{G\times G}(D_G(Y)\wedge
                                                     Y)}\longrightarrow 
    {(\on{EG}\times\on{EG})\ltimes_{G}(D_G(Y)\wedge Y)}.
  $$
  The second map is $$(\on{EG}\times\on{EG})\ltimes_G-$$ applied to
  the canonical $G$-map 
  \begin{equation}
    \label{auswertung-Eqn}
    \beta_G \negmedspace : D_GY\wedge Y \longrightarrow \SS
  \end{equation}
  (cf.\ Theorem \ref{strong-dual-Thm}).
  Its target is (the suspension spectrum of)
  $$
    (\on{EG}\times\on{EG})_+\wedge_{G}\SS\simeq\BG_+.
  $$
  The last map is 
  $$
    (Bp_G)_+ : \BG_+\longrightarrow \SS,
  $$
  where $p_G$ denotes the unique map from $G$ to the trivial group.
  This completes the construction of the map (\ref{adjoint-Eqn}).
%
\begin{Rem}
  \label{sphere-Obs}
  In the case $Y=\SS$, the construction of (\ref{adjoint-Eqn})
  specializes to the definition of Strickland's
  inner product \cite[8.2]{Strickland:inner:product}
  $$
    \longmap{\beta}{\BG_+\wedge\BG_+}{\SS}.
  $$ 
\end{Rem}
\begin{Cor}
  Theorem \ref{Borel-Thm} is true for $Y=\SS$.
\end{Cor}
\begin{Pf}{}
  This is the fact that the Strickland inner product is non-degenerate
  \cite[8.3]{Strickland:inner:product}.
\end{Pf}
The second easiest special case of Theorem \ref{Borel-Thm} is the case
that $Y$ is a different zero sphere. 
\begin{Prop}\label{zero-sphere-Prop}
  For $Y=\zs$ the map (\ref{adjoint-Eqn}) is the Strickland inner
  product on $$(\borel G/H)_+,$$ and the theorem holds for $Y=\zs$.
\end{Prop} 
\begin{Pf}{}
  Recall from \cite[p.176]{May:CBMS} that for finite groups $H\sub G$
  a $G$-equivariant (strong) dual of $\zs$ is $\zs$, with the map
  $\beta_G$ in 
  (\ref{auswertung-Eqn}) given by the composite (of space level maps)
  \begin{equation}\label{G/H-Eqn}
    (G/H\times G/H)_+\longto\zs\longto\SS,
  \end{equation} 
  where the first map is the $G$-equivariant Pontrjagin-Thom collapse along the
  diagonal inclusion (that is, it is the identity on the diagonal and
  everything else gets 
  mapped to the basepoint) and the second map is $p_+$, where $p$ is
  the unique ($G$-equivariant) map 
  $$\longmap p{G/H}\pt.$$
  Following Strickland, we write $\BGH$ for the Borel construction
  of the finite groupoid $\mathcal G$ defined by the action of $G$ on $G/H$,
  remembering that
  $$
    (\BGH)_+ = \Borel \zs.
  $$
  Strickland defines the inner product on (the suspension spectrum of)
  $(\BGH)_+$ as the composite
  $$
    \beta\negmedspace 
    : (B\mathcal G\times B\mathcal G)_+
    \stackrel{T\delta}{\longrightarrow} (\BGH)_+
    {\longrightarrow}\SS,   
  $$
  where 
  $$
    \delta = \delta_{\gpd}\negmedspace 
    : \mathcal G \longto 
    \mathcal G \times \mathcal G
  $$
  is the diagonal inclusion of groupoids, and the second map is the
  Borel construction of the unique map of groupoids
  $$
    \longmap{p_{\gpd}}{\gpd}{1}.
  $$
  Note that $\delta_{\gpd}$ factors as
  $$
    \delta_{\gpd}\negmedspace : (G/H\acted G)
    \stackrel{i}{\hookrightarrow} (G/H\times G/H)\acted G
    {\hookrightarrow} 
    (G/H\times G/H)\acted G\times G, 
  $$
  where $i$ is an inclusion of finite $G$-sets (namely the diagonal
  inclusion mentioned above), and the second map is the diagonal
  inclusion of groups $\delta_G$, whose
  transfer $T\delta_G$ is the
  first map in the construction of (\ref{adjoint-Eqn}). 
  We need to identify
  $\on{Ti}$. However, $\on{Bi}$ is a particularly simple example of
  covering with 
  finite fibers, namely the inclusion of some path components. A look
  at the construction of $\on{Ti}$ in \cite[4.1.1]{Adams:infinite} shows
  that $\on{Ti}$ is given by the (space level) map
  $$
    (\borel(G/H\times G/H))_+\longrightarrow (\borel G/H)_+
  $$
  that is the identity on $\on{im}(\on{Bi})$ and maps everything else to
  the basepoint. This is exactly $\Borel-$ applied to the
  Pontrjagin-Thom collapse in (\ref{G/H-Eqn}). The map $p_{\gpd}$ 
  factors as 
  $$
    p_{\gpd}\negmedspace : (G/H\acted G)\stackrel{p}{\longto}(\pt\acted G)
    \stackrel{p_G}{\longrightarrow} (\pt\acted 1),
  $$
  where $p$ is as in (\ref{G/H-Eqn}). Together this proves the claim
  that (\ref{adjoint-Eqn}) is the Strickland inner product:
  $$
    \beta = (\on{Bp_{\gpd}})_+\circ\on{T\delta_{\gpd}} =
    (\on{Bp_G})_+\circ(\on{Bp})_+\circ\on{Ti}\circ\on{T\delta_G} = 
    (\on{Bp_G})_+\circ(\Borel\beta_{G})\circ\on{T\delta_G}. 
  $$
  As above, non-degeneracy of the Strickland product in the
  $\Kn$-local category implies that (\ref{Borel-Eqn}) is an
  isomorphism for $Y=\zs$.  
\end{Pf}
Before we proceed to higher dimensional spheres, we recall that in any
closed symmetric monoidal category we have an isomorphism
\begin{equation}
  \label{smash-D-Eqn}
  D(X)\wedge D(Y) \stackrel\cong\longrightarrow D(X\wedge Y), 
\end{equation}
which identifies the evaluation map $\beta_{X\wedge Y}$ with
$$
  (\beta_X\wedge\beta_Y)\circ(\id_{DX}\wedge \tau\wedge\id_{Y}), 
$$
where $\tau$ switches $DY$ and $X$.
\begin{Lem}\label{sphere-Lem}
  Theorem \ref{Borel-Thm} holds for spheres
  $$
    Y = \zs\wedge\SSn.
  $$
\end{Lem}
\begin{Pf}{}
  The sphere $\mathbb{S}^n$ is strongly dualizable in $\mathcal S$ with
  dual $\SSmn$, and both functors $\mathcal S\to\mathcal S_G$
  and $\mathcal S\to\mathcal S_\Kn$ preserve the data of strong
  dualizability in Theorem \ref{strong-dual-Thm}. By
  (\ref{smash-D-Eqn}) we have
  $$
    D_G(\zs\wedge\SSn)\cong D_G(\zs)\wedge\SSmn.
  $$
  Since $\SSn$ and $\SSmn$ have trivial $G$-action, we have
  $$
    \thsG(\zs\wedge\mathbb{S}^{\pm n}) = (\Borel\zs)\wedge\mathbb{S}^{\pm n}, 
  $$
  and under this identification the map (\ref{adjoint-Eqn}) becomes
  $$
    \beta_{\Borel(\zs)}\wedge\beta_\SSn :
    \Borel(D_G(\zs))\wedge\Borel(\zs)\wedge \SSn\wedge\SSmn\longto\SS. 
  $$
  Here $\beta_{\Borel(\zs)}$ is the map of the theorem for $\zs$, and
  by Proposition \ref{zero-sphere-Prop}, it is the evaluation map
  of a strong duality in $\mathcal S_\Kn$. The map $\beta_\SSn$ is
  already an evaluation of a strong duality in $\mathcal S$ and thus
  also in $\mathcal S_\Kn$. We apply (\ref{smash-D-Eqn}) again, this time
  in the $\Kn$-local category, to complete the
  proof.   
\end{Pf}

\begin{Pf}{of Theorem \ref{Borel-Thm}}
  We prove the theorem by induction over the cells.
  All our categories have compatible triangulated and closed symmetric
  monoidal structures. In particular, duals commute with direct sums and
  take triangles (in the opposite category) into triangles.
  The Borel construction also preserves the triangulated structure.
  Since both sides of (\ref{Borel-Eqn}) preserve finite sums,
  Lemma \ref{sphere-Lem} implies that the statement is true for
  finite bouquets of spheres.
  Both sides of (\ref{Borel-Eqn}) preserve triangles, thus if the
  theorem is true for two objects in an exact triangle, it is also
  true for the third.
\end{Pf}

\begin{Pf}{{of Theorem \ref{dual-independent-Thm}}}
The map of Theorem \ref{dual-independent-Thm} is the Borel construction
$$\on{Borel}(p_{M\mmod G})_+$$
of the unique map of orbifolds from $M\mmod G$ to a point. This map
factors as
$${p_{M\mmod G}}:{M\mmod G}\stackrel{\pi}{\longrightarrow} {\pt\mmod
  G}\stackrel{p_G}{\longrightarrow}\pt,$$
where $\pi$ is the unique $G$-map from $M$ to a point, and $p_G$ is
the unique map from $G$ to the trivial group.
Recall from Theorem \ref{Neil-Thm} that $\eta$ is defined as 
$$
  \eta\negmedspace :\SS\stackrel\cong\longrightarrow D\SS {\longrightarrow}
  D(\BG_+)\stackrel\cong\longleftarrow\BG_+, 
$$
where the second map is $D((Bp_G)_+)$, and
the last map is the adjoint $\beta^\sharp$ of the Strickland inner
product on $\BG_+$.

By \cite[XVI.8.1]{May:CBMS}, the Thom spectrum $M^{-\tau}$ is a
$G$-equivariant (strong) dual of $M_+$, and the $G$-equivariant \PT
collapse 
$$
  \longmap{PT}{\SS\acted G}{M^{-\tau}\acted G}
$$
is the dual  
$D_G({\pi_+})$. 
Under the isomorphism of Theorem \ref{Borel-Thm} (vertical arrows in
the diagram below), the
twisted half smash product $\thsG(PT)$ (top row) becomes     
$$
  \xymatrix{{\thsG (D_GM_+)}\ar_\cong[d]& {\Borel (D_G\SS)}\ar^\cong[d]\ar[l]&
  {\Borel\SS}\ar_{\phantom{Xx}\cong}[l]\\
  {D(\Borel M_+)}& D(\Borel\SS),\ar[l]}
$$
where the bottom arrow is $D(\Borel\pi_+)$.
By Remark \ref{sphere-Obs}, the composite of the two rightmost arrows
is $\beta^\sharp$.
When precomposing with $\eta$, $\beta^\sharp$ and its inverse cancel out,
and we obtain
$$
  \thsG(PT)\circ\eta = D((\borel\pi)_+)\circ D((Bp_G)_+) =
  D(\on{Borel}(p_{M\mmod G})_+),
$$
which completes the proof.
\end{Pf}
%
%
%
%
\section{The DMVV formula}
\subsection{Conjugacy classes of $h$-tuples of commuting elements of $\Sl$}
\label{h-tup-Sec}
Just as the conjugacy classes of elements of $\Sl$ are in one to one
correspondence with partitions
$$\sum a_nn = l$$
(i.e. the shape of the Young tableau), one also describes conjugacy
classes of $h$-tuples of commuting elements in terms of the
corresponding orbit decomposition of the set
$$\ll:=\{1,\dots,l\}.$$
More precisely, such an $h$-tuple $\gg$ defines an action of $\ZZ^h$
on $\ll$, and $\ll$ decomposes into orbits of that action. Two such
$h$-tuples are conjugate by a permutation $g$ of $\ll$ if and only if
their orbit decompositions are isomorphic (and an isomorphism is given
by $g$).

Orbits are finite transitive $\ZZ^h$-sets, and every finite transitive
$\ZZ^h$-set $\T$ turns up as a possible orbit for $l\geq
|\T|$, where $|\T|$ is the number of elements in $\T$. 

Let $\mathcal{T}=\{\T\}$ contain one representative for each
isomorphism class of finite transitive $\ZZ^h$-sets. The above
discussion summarizes as follows. The conjugacy classes of $h$-tuples
of commuting elements in $\Sl$ are classified by expressions
$$\sum_{\TT}a_\T\T\text{  s.t.  }\sum_\TT a_\T|\T| = l,$$
where for given $\gg$ the expression
$\sum_{\TT}a_\T\T$ counts the number $a_\T$ of times each isomorphism
class of finite transitive $\ZZ^h$-set $\T$ occurs in the 
decomposition of $\ll$ into orbits of the subgroup $\langle
g_1,\dots,g_h\rangle$ generated by the $g_i$.
If the conjugacy class $[g_1,\dots,g_h]$ corresponds to $\sum_\TT
a_\T\T$, then  
the centralizer of $\alpha$ in $G$ can be described as follows
\begin{equation}\label{centralizer}
  C_\gg\cong\prod_\TT\Aut_{\Zh}(\T)^{a_\T}\rtimes\Sigma_{a_\T},
\end{equation}
where $\Sigma_{a_\T}$ permutes the $a_\T$ orbits isomorphic to $\T$ and 
$\Aut_\Zh(\T)$ acts on each of them individually.
Now $\T$ is a transitive $\Zh$-set, and $\Zh$ is abelian. This means
that for any $x,y\in\T$ there is an element $z\in\Zh$ such that $zx =
y$, and multiplication with $z$ is the unique automorphism of $\T$
mapping $x$ to $y$. Thus the number of elements in $\Aut_\Zh(\T)$ is
$$
  |\Aut_\Zh(\T)| = |\T|.
$$
Since the conjugacy class of $\gg$ in $\Hom(\Zh,\Sl)$ is the orbit of
$\gg$ under the action of $\Sl$ by conjugation, we have
$$
  [g_1,\dots,g_h]_\Sl\cong_\Sl\Sl/C_\gg.
$$
Therefore, its number of elements is by (\ref{centralizer})
$$
  |[g_1,\dots,g_h]_\Sl| = \frac{l!}{\prod_\TT|\T|^{a_\T}a_\T}.
$$
Assume now that we are only interested in $h$-tuples of commuting
elements of $p$-power order. Then the same discussion goes through,
but we need to replace $\mathcal T$ by
the set $\mathcal T_p$ containing one representative for each
isomorphism class of finite transitive $\ZZ_p^h$-set. Note that
elements of $\mathcal T_p$ have $p$-power cardinalities, since each of
them
can be identified with a quotient of $(\ZZ/p^j\ZZ)^h$ for some
sufficiently large $j$.
\subsection{The DMVV formula for $\phi_\orb$}
\label{proof-of-main-thm-Sec} 
We start by proving a formula for the total symmetric power.
Let $S_t$ be as in Definition \ref{higher-symmetric-powers-Def} and
$T_{p^k}$ as in Definition \ref{Hecke-Def}.
\begin{Prop}\label{higherchromaticEquation}
  We have
  $$
    S_t(x) = \exp \left[ \sum_{k\geq 0} T_{p^k}(x)t^{p^k} \right].
  $$
\end{Prop}
\begin{Pf}{}
We have
$$\exp \left[ \sum_{k\geq 0} {T_{p^k}(x)}t^{p^k} \right] = 
\sum_{m\geq 0}\frac{1}{m!}
\left[ \sum_{k\geq 0} {T_{p^k}(x)}t^{p^k} \right]^m.$$
In this equation, the coefficient of $t^l$ is
$$\sum_{l = \sum\limits_{\TTp} a_\T|\T|}\frac{(\sum a_\T)!}{\prod(a_\T!)}
\frac{1}{(\sum a_\T)!}\prod_{\TTp}\left(\frac{\psi_\T(x)}{|\T|}\right)^{a_\T},$$
where $\frac{(\sum a_\T)!}{\prod(a_\T!)}$ counts the number of ways to
partition a 
set of $\sum a_\T$ (orbits) into subsets of size $a_\T$ (the number
of times $\T$ occurs as orbit), and $\frac{1}{(\sum a_\T)!}$ is
$\frac{1}{m!}$. This is
\begin{eqnarray*}
\sum_{\sum a_\T|\T|=l}\prod_{\mathcal{T}_p}
\frac{1}{|\T|^{a_\T}(a_\T!)}\psi_\T(x)^{a_\T} 
& = & \sum_{\sum a_\T|\T|=l}\left( \prod_{\mathcal{T}_p}
\frac{1}{|\T|^{a_\T}(a_\T!)}\right)
 \psi_{(\coprod\limits_{\mathcal{T}_p} a_\T\T)}(x) \\
& = & \sum_{[\alpha]}\GG{C_\alpha}\psi_\alpha(x)\\
& = & \sigma_l(x).
\end{eqnarray*}
\end{Pf}
We are now ready to prove Theorem \ref{main-Thm} of the introduction.
\begin{Thm}\label{DMVV-Thm}
  Let $\phi$ be an $\hinf$-orientation of $E_h$. Then
  $$
    \sum_{n\geq 0}\phi_{\on{orb}}(M^n\mmod\Sn) t^n
    = \exp \left[ \sum_{k\geq 0} T_{p^k}(\phi(M))t^{p^k} \right].
  $$
\end{Thm}
\begin{Pf}{}
  We have 
  \begin{eqnarray*}
    \sum_{n\geq 0}\phi_{\on{orb}}(M^n\mmod\Sn) t^n &=&
      \sum_{n\geq 0}\eta_\Sn^*\circ\phi_{\Sn}(M^n) t^n \\
    & = & \sum_{n\geq 0}\eta_\Sn^*\circ \phi_{\Sn}\circ P_n^{\MU}(M)t^n \\
    & = & \sum_{n\geq 0}\eta_\Sn^*\circ P_n^{E_h}\circ\phi(M)t^n\\
    & = & S_t(\phi(M)),
  \end{eqnarray*}
  where the first equation is Definition \ref{orb-Def}, $P_n^{\MU}$ and
  $P_n^{E_h}$ are the $n^{th}$ power operations in cobordism and Morava
  $E$-theory of a point, the third equation holds,
  because $\phi$ is an $\hinf$-map, and the fourth equation is the
  definition of $S_t$. The claim now follows from Proposition
  \ref{higherchromaticEquation}. 
\end{Pf}
Note the striking similarity of the right hand side of the
DMVV-formula with the formal inverse of Rezk's logarithm formula
\cite[p.4]{Rezk:logarithms}
$$
\exp\sum_{k\geq 0} T_{p^k}(-).
$$
Here the $T_{p^k}$ are as in Definition \ref{Hecke-Def}.
%
%
%
%
\begin{Exa}[{$\sigma$-orientation}]
  Any elliptic spectrum $E$ has a canonical orientation
  $$\longmap{\sigma_E}{\MU\langle6\rangle}{E},$$
  and it was shown in \cite{AHSII}, that in the case $E=E_2$, the map
  $\sigma$ is an $\hinf$-map.
\end{Exa}
  The following result due to Ando classifies the complex genera into
  $E_h$ that can be taken as input for Theorem \ref{main-Thm}:
\begin{Thm}[{\cite{Ando:Duke}}]
  The spectrum $E_h$ is an $\hinf$-spectrum. A map of ring spectra
  $$\longmap\phi\MU{E_h}$$
  is an $\hinf$-map if and only if the $p$-series of its Euler class
  $e_\phi$ (of the universal line bundle) satisfies 
  $$
    \left[p\right]_F(e_\phi) = \prod_{\stackrel{v\in F(D_1)}{[p]_F(v)=0}}(v+_Fe_\phi),
  $$
  where $D_1$ denotes the ring
  extension of $E^0_n$ obtained by adjoining 
  the roots of the $p$-series of $F$, and $F(D_1)$ stands for the
  maximal ideal of $D_1$ with the group structure $x+_Fy$.
\end{Thm}
\subsection{Atiyah-Tall-Grothendieck type definition of Hecke
  operators}\label{Atiyah-Tall-Sec} 
  The left-hand side of the equation in  Proposition
  \ref{higherchromaticEquation} is defined in greater generality than
  its right-hand side, motivating the following definition.
  Let $E$ be an even periodic, $\Kn$-local $\hinf$-spectrum.
  Then the total symmetric power $S_t$ is defined on elements of
  $E(X)$ and takes values in
  $$
    1+tE(X)\ps{t}
  $$
  (cf.\ Definition \ref{generalized-symmetric-powers-Def}).
  \begin{Def}
    In this situation we define additive operators $T_n$ on $E(X)$ by
    $$
      \sum_{n\geq 1}T_nt^n := \log S_t.
    $$    
  \end{Def}
  Following Grothendieck \cite{Grothendieck}, or the interpretation
  for $K$-theory by Atiyah and Tall \cite{Atiyah:Tall}, we note that 
  $$
    t\frac{d}{\on{dt}}\log S_t(x) =
    t\frac{\frac{d}{\on{dt}}S_t(x)}{S_t(x)}     
  $$
  takes values in $E(X)\ps t$.
  Thus the Hecke operators are operations
  $$\longmap{T_n}{E(X)}{\frac{1}{n}E(X)}.$$
  We can make the connection to the Atiyah-Tall-Grothendieck
  definition of the Adams operations even more precise:  
  Let
  $$
    \Lambda_{t} := \frac{1}{S_{-t}}
  $$
  denote the ``total exterior power'' in $E$-theory. 
  This defines a $\lambda$-ring structure on $E(X)$, whose Adams
  operations are given by $\psi_n = nT_n$. 
%
%
%
\bibliographystyle{alpha}
\bibliography{thesis}
\end{document}